\documentclass[twoside,10pt]{amsart}
\usepackage{amssymb}
\usepackage{graphicx}
\usepackage{amsmath}
\usepackage{amsfonts}

\usepackage[latin1]{inputenc}
\usepackage[T1]{fontenc}
\usepackage[english]{babel}

\newcommand{\C}{\mathbb{C}}
\newcommand{\R}{\mathbb{R}}
\newcommand{\K}{\mathbb{K}}
\newcommand{\N}{\mathbb{N}}

\newcommand{\T}{\mathbb{T}}
\newcommand{\eps}{\varepsilon}
\newcommand{\ext}[1][X^*]{\ensuremath{\mathrm{ext}(B_{#1})}}
\newcommand{\extr}{\ensuremath{\mathrm{ext}}}
\newcommand{\e}{\mathrm{e}}
\newcommand{\Id}{\mathrm{Id}}

\newtheorem{theorem}{\sffamily Theorem}
\newtheorem{proposition}{\sffamily Proposition}
\newtheorem{lemma}{\sffamily Lemma}
\newtheorem{corollary}{\sffamily Corollary}
\newtheorem{example}{\sffamily Example}

\newtheorem{problem}{\sffamily Problem}
\newtheorem{examples}[example]{\sffamily Examples}

\def\qed{\quad{$\blacksquare$}\medskip }

 \DeclareMathOperator{\re}{Re}
 \DeclareMathOperator{\dist}{dist}
 \DeclareMathOperator{\ecc}{\overline{\mathrm{co}}}
 \DeclareMathOperator{\ec}{co}

\begin{document}

\title[Numerical index of Banach spaces]{Recent progress and open questions on
 \\ the numerical index of Banach spaces}

\date{March 27th, 2006; minor modifications May 24th, 2006}

\author[V.~Kadets]{Vladimir Kadets}
\address[Vladimir Kadets]{Faculty of Mechanics and Mathematics\\
Kharkov National University\\ pl.\ Svobody 4\\ 61077 Kharkov,
Ukraine \newline {\tt vova1kadets@yahoo.com}}

\author[M.~Mart\'{\i}n]{Miguel Mart\'{\i}n}
\address[Miguel Mart\'{\i}n \and Rafael Pay\'{a}]{Departamento de An\'{a}lisis
Matem\'{a}tico\\ Facultad de Ciencias\\ Universidad de Granada\\ 18071
Granada, Spain \newline {\tt mmartins@ugr.es, \ \ \ rpaya@ugr.es}}

\author[R.~Pay\'{a}]{Rafael Pay\'{a}}

\begin{abstract}
The aim of this paper is to review the state-of-the-art of recent
research concerning the numerical index of Banach spaces, by
presenting some of the results found in the last years and proposing
a number of related open problems.
\end{abstract}
\maketitle

\thispagestyle{empty}

\section*{Introduction}
The numerical index of a Banach space is a constant relating the
behavior of the numerical range with that of the usual norm on the
Banach algebra of all bounded linear operators on the space. The
notion of numerical range (also called field of values) was first
introduced by O.~Toeplitz in 1918 \cite{Toe} for matrices, but his
definition applies equally well to operators on infinite-dimensional
Hilbert spaces. If $H$ denotes a Hilbert space with inner product
$(\cdot\mid\cdot)$, the \emph{numerical range} of a bounded linear
operator $T$ on $H$ is the subset $W(T)$ of the scalar field defined
by
$$
W(T):=\big\{(T x\mid x) \ : \ x\in H,\ (x\mid x)=1\big\}.
$$
Some properties of the Hilbert space numerical range are discussed
in the classical book of P.~Halmos \cite[\S17]{Hal}. Let us just
mention that the numerical range of a bounded linear operator is
(surprisingly) convex and, in the complex case, its closure contains
the spectrum of the operator. Moreover, if the operator is normal,
then the closure of its numerical range coincides with the convex
hull of its spectrum. Further developments can be found in a recent
book of K.~Gustafson and D.~Rao \cite{GusRao}. In the sixties, the
concept of numerical range was extended to operators on general
Banach spaces by G.~Lumer \cite{Lumer} and F.~Bauer \cite{Bauer}.
Let us give the necessary definitions. Given a real or complex
Banach space $X$, we write $B_X$ for the closed unit ball and $S_X$
for the unit sphere of $X$. The dual space will be denoted by $X^*$
and $L(X)$ will be the Banach algebra of all bounded linear
operators on $X$. The \emph{numerical range} of an operator $T\in
L(X)$ is the subset $V(T)$ of the scalar field defined by
$$
V(T):=\{x^*(Tx)\ : \ x\in S_X,\ x^*\in S_{X^*},\ x^*(x)=1\}.
$$
The \emph{numerical radius} is the seminorm defined on $L(X)$ by
$$
v(T):=\sup\{|\lambda|\ : \ \lambda\in V(T)\}
$$
for $T\in L(X)$. Classical references here are the monographs by
F.~Bonsall and J.~Duncan \cite{B-D1,B-D2} from the seventies. Let us
mention that the numerical range of a bounded linear operator is
connected (but not necessarily convex, see
\cite[Example~21.6]{B-D2}) and, in the complex case, its closure
contains the spectrum of the operator. The theory of numerical
ranges has played a crucial role in the study of some algebraic
structures, especially in the non-associative context (see the
expository paper \cite{K-M-RP-survey} by A.~Kaidi, A.~Morales, and
A.~Rodr\'{\i}guez Palacios, for example).

The concept of \emph{numerical index} of a Banach space $X$ was
first suggested by G.~Lumer in 1968 (see \cite{D-Mc-P-W}), and it is
the constant $n(X)$ defined by
$$
n(X):=\inf\{v(T)\ : \ T\in L(X),\ \|T\|=1\}
$$
or, equivalently,
$$
n(X)=\max\{k\geqslant 0 : \ k\,\|T\|\leqslant v(T) \ \ \forall\,
T\in L(X)\}.
$$
Note that $n(X)>0$ if and only if $v$ and $\|\cdot\|$ are equivalent
norms on $L(X)$. At that time, it was known that in a complex
Hilbert space $H$ with dimension greater than $1$, $\|T\|\leqslant 2
v(T)$ for all $T\in L(H)$, and $2$ is the best constant; in the real
case, there exists a norm-one operator whose numerical range reduces
to zero. In our terminology, $n(H)=1/2$ if $H$ is complex, and
$n(H)=0$ if it is real. Actually, real and complex general Banach
spaces behave in a very different way with regard to the numerical
index, as summarized in the following equalities by J.~Duncan,
C.~McGregor, J.~Pryce, and A.~White \cite{D-Mc-P-W}:
\begin{eqnarray}\label{eq-DMCPW}
\{n(X)\ : \ X \text{ complex Banach space} \ \}&=&[\e^{-1},1], \\
\{n(X)\ : \ X \text{ real Banach space} \ \}&=&[0,1].\notag
\end{eqnarray}
The fact that $n(X)\geqslant \e^{-1}$ for every complex Banach space
$X$ was observed by B.~Glickfeld \cite{Gli} (by making use of a
classical theorem of H.~Bohnenblust and S.~Karlin \cite{B-K}), who
also gave an example where this inequality becomes an equality. It
is showed in the already cited paper \cite{D-Mc-P-W}, that
$M$-spaces, $L$-spaces and their isometric preduals, have numerical
index~$1$, a property shared by the disk algebra (M.~Crab,
J.~Duncan, and C.~McGregor \cite[Theorem~3.3]{C-D-Mc}). Finally, let
us mention that the real space $X_\R$ underlying a complex Banach
space $X$ satisfies $n(X_\R)=0$; actually, the isometry
$x\longmapsto i\,x$ has numerical radius $0$ when viewed as an
operator on $X_\R$.

In the last ten years, many results on the numerical index of Banach
spaces have appeared in the literature. This paper aims at reviewing
the state of the art on this topic and proposing a variety of open
questions. The structure of our discussion is as follows.

\tableofcontents

We finish this introduction by recalling some definitions and fixing
notation. We write $\T$ to denote the unit sphere of the base field
$\K$ ($=\R$ or $\C$). We use the notation $\re(\,\cdot\,)$ to denote
the real part function, which should be considered as the identity
when $\K=\R$. Given a real or complex Banach space $X$, we write
$\ec(B)$ and $\ecc(B)$ to denote, respectively, the convex hull and
the closed convex hull of a set $B\subseteq X$ and we denote by
$\extr(C)$ the set of extreme points of a convex set $C\subseteq X$.
A subset $A$ of $B_{X^*}$ is said to be \emph{norming (for $X$)} if
$$
\|x\|=\sup\{|x^*(x)| \ : \ x^*\in A\} \qquad \big(x\in X\big)
$$
or, equivalently, if $B_{X^*}=\ecc^{w^*}(\T\,A)$ (Hahn-Banach
Theorem). Finally, for $1\leqslant p \leqslant\infty$, we write
$\ell_p^m$ to denote the normed space $\K^m$ endowed with the usual
$p$-norm, and we write $X \oplus_p Y$ to denote the $\ell_p$-direct
sum of the spaces $X$ and $Y$.

\section{Computing the numerical index}\label{sec:computing}
In view of the examples given in the introduction, the most
important family of classical Banach spaces (in the sense of
H.~Lacey \cite{Lac}) whose numerical indices remain unknown is the
family of $L_p$ spaces when $p\neq 1,2,\infty$. This is actually one
of the most intriguing open problems in the field but, very
recently, E.~Ed-Dari and M.~Khamsi \cite{Ed-Dari,Ed-Dari-Khamsi}
have made some progress. We summarize their results in the following
statement and  use it to motivate some conjectures.

\begin{theorem}[\mbox{\rm \cite{Ed-Dari,Ed-Dari-Khamsi}}] \label{th:lp}
Let $1\leqslant p\leqslant \infty$ be fixed. Then,
\begin{enumerate}
\item[(a)] $n(L_p[0,1])=n(\ell_p)=\inf\{n(\ell_p^m)\, : \, m\in
\N\}$, and the sequence $\{n(\ell_p^m)\}_{m\in\N}$ is decreasing.
\item[(b)] $n(L_p(\mu))\geqslant n(\ell_p)$ for every positive measure $\mu$.
\item[(c)] In the real case,
$$
\frac12 M_p \leqslant n(\ell_p^2) \leqslant M_p, \qquad \text{where
 }\ M_p=\sup_{t\in[0,1]}\frac{|t^{p-1}-t|}{1+t^p}.
$$
\end{enumerate}
\end{theorem}

When $p\neq 2$, it is known that $n(\ell_p^m)>0$ for every
$m\geqslant 2$, and also that $v(T)>0$ for every non-null $T\in
L(\ell_p)$ (see \cite[Theorem~2.3 and subsequent
remark]{Ed-Dari-Khamsi}), but we do not know if $n(\ell_p)>0$.
Observe that a positive answer to this question implies, thanks to
(b) above, that $n(L_p(\mu))>0$ for every positive measure $\mu$.

\begin{problem}
Is $n(\ell_p)$ positive for every $p\neq 2$?
\end{problem}

With respect to item (c) in the above theorem, let us explain the
meaning of the number $M_p$. It can be deduced from
\cite[\S3]{D-Mc-P-W} that, given an operator $T\in L(\ell_p^2)$
represented by the matrix $\begin{pmatrix} a & b
\\ c & d \end{pmatrix}$, one has
\begin{align}\label{eq:numradlp}
v(T)=\max &\left\{
\underset{z\in\T}{\underset{t\in[0,1]}{\max}}\frac{\left|a+d\,t^p+z\,b\,t+
\overline{z}\,c\,t^{p-1}\right|}{1+t^p},\
\underset{z\in\T}{\underset{t\in[0,1]}{\max}}
\frac{\left|d+a\,t^p+\overline{z}\,c\,t+z\,b\,t^{p-1}\right|}{1+t^p}
\right\}.
\end{align}
It follows that $M_p$ is equal to the numerical radius of the
norm-one operator $U\equiv \begin{pmatrix} 0 & 1
\\ -1 & 0 \end{pmatrix}$ in $L(\ell_p^2)$ (real case), so $n(\ell_p^2)\leqslant M_p$.
For $p=2$, the operator $U$ has minimum numerical radius, namely
$0$. We may ask if $U$ is also the norm-one operator with minimum
numerical radius for all the real spaces $\ell_p^2$.

\begin{problem}
Is it true that, in the real case, $\displaystyle
n(\ell_p^2)=\sup_{t\in[0,1]}\frac{|t^{p-1}-t|}{1+t^p}$ for every $1<
p < \infty$?
\end{problem}

In the complex case, the operator $U$ acting on $\ell_2^2$ satisfies
$v(U)=\|U\|$ (take $z=i$ and  $t=1$ in \eqref{eq:numradlp}) and,
therefore, its numerical radius is not the minimum. Actually, one
has
$$
n(\ell_2^2)=\frac12=v(S),
$$
where $S\in L(\ell_2^2)$ is the `shift' $S\equiv
\begin{pmatrix} 0 & 0
\\ 1 & 0 \end{pmatrix}$. Therefore, we bet that $n(\ell_p^2)=v(S)$
for every $p$ in the complex case. It can be checked from
\eqref{eq:numradlp} that
$$
v(S)=\dfrac{(p-1)^{\frac{p-1}{p}}}{p}=\dfrac{1}{p^{\frac1p}\,q^{\frac1q}},
$$
where $\frac1p +\frac1q=1$.

\begin{problem}
Is it true that, in the complex case,
$n(\ell_p^2)=\dfrac{1}{p^{\frac1p}\,q^{\frac1q}}$ for every
$1<p<\infty$?
\end{problem}

In view of Theorem~\ref{th:lp}.a, the two-dimensional case is only
the first step in the computation of $n(\ell_p)$, but it is
reasonable to expect that the sequence $\{n(\ell_p^m)\}_{m\in \N}$
is always constant, as it happens in the cases $p=1,2,\infty$.

\begin{problem}
Is it true that $n(\ell_p)=n(\ell_p^2)$ for every $1<p<\infty$?
\end{problem}

In a 1977 paper \cite{Hur}, T.~Huruya determined the numerical index
of a $C^*$-algebra. Part of the proof was recently clarified by
A.~Kaidi, A.~Morales, and A.~Rodr\'{\i}guez-Palacios in \cite{K-M-RP},
where the result is extended to $JB^*$-algebras and preduals of
$JBW^*$-algebras. Let us state here those results just for
$C^*$-algebras and preduals of von Neumann algebras.

\begin{theorem}[\mbox{\rm \cite{Hur} and
\cite[Proposition~2.8]{K-M-RP}}] \label{th:cstaralgebras} Let $A$ be
a $C^*$-algebra. Then, $n(A)$ is equal to $1$ or $\frac12$ depending
on whether or not $A$ is commutative. If $A$ is actually a von
Neumann algebra with predual $A_*$, then $n(A_*)=n(A)$.
\end{theorem}

We do not know if there is an analogous result in the real case. We
recall that a \emph{real $C^*$-algebra} can be defined as a
norm-closed self-adjoint real subalgebra of a complex $C^*$-algebra,
and a \emph{real $W^*$-algebra} (or \emph{real von Neumann algebra})
is a real $C^*$-algebra which admits a predual (see \cite{IsiRod}
for more information).

\begin{problem}
Compute the numerical index of real $C^*$-algebras and isometric
preduals of real $W^*$-algebras.
\end{problem}

The fact that the disk algebra has numerical index~$1$ was extended
to function algebras by D.~Werner in 1997 \cite{WerJFA}. A
\emph{function algebra} $A$ on a compact Hausdorff space $K$ is a
closed subalgebra of $C(K)$ which separates the points of $K$ and
contains the constant functions.

\begin{proposition}[\mbox{\rm \cite[Corollary~2.2 and proof of
Theorem~3.3]{WerJFA}}] If $A$ is a function algebra, then $n(A)=1$.
\end{proposition}

Of course, there are many other Banach spaces whose numerical index
is unknown. We propose to calculate some of them.

\begin{problem}
Compute the numerical index of $C^m[0,1]$ (the space of $m$-times
continuously differentiable real functions on $[0,1]$, endowed with
any of its usual norms), $\text{Lip}(K)$ (the space of all Lipschitz
functions on the complete metric space $K$), Lorentz spaces, and
Orlicz spaces.
\end{problem}

Some of the classical results given in the introduction about the
numerical index of particular spaces have been extended to sums of
families of Banach spaces and to spaces of vector-valued functions
in various papers by G.~L\'{o}pez, M.~Mart\'{\i}n, J.~Mer\'{\i}, R.~Pay\'{a}, and
A.~Villena \cite{L-M-M,M-P,M-V}.

We start by presenting the result for sums of spaces. Given an
arbitrary family $\{X_\lambda\, : \, \lambda\in\Lambda\}$ of Banach
spaces, we denote by $\left[\oplus_{\lambda\in\Lambda}
X_\lambda\right]_{c_0}$, $\left[\oplus_{\lambda\in\Lambda}
X_\lambda\right]_{\ell_1}$ and $\left[\oplus_{\lambda\in\Lambda}
X_\lambda\right]_{\ell_\infty}$ the $c_0$-, $\ell_1$- and
$\ell_\infty$-sum of the family.

\begin{proposition}[\mbox{\rm \cite[Proposition~1]{M-P}}] \label{suma}
Let $\{X_\lambda:\ \lambda \in \Lambda\}$ be a family of Banach
spaces. Then
\begin{equation*}
n\Bigl(\left[\oplus_{\lambda\in\Lambda} X_\lambda\right]_{c_0}
\Bigr)= n\Bigl(\left[\oplus_{\lambda\in\Lambda}
X_\lambda\right]_{\ell_1} \Bigr)=
n\Bigl(\left[\oplus_{\lambda\in\Lambda}
X_\lambda\right]_{\ell_\infty} \Bigr)= \inf_{\lambda}
\,n(X_\lambda).
\end{equation*}
\end{proposition}

The above result is not true for $\ell_p$-sums if $p$ is different
from $1$ and $\infty$. Nevertheless, it is possible to give one
inequality and, actually, the same is true for absolute sums
\cite{MarEMA}. Recall that a direct sum $Y\oplus Z$ is said to be an
\emph{absolute sum} if $\|y + z\|$ only depends on $\|y\|$ and
$\|z\|$ for $(y,z)\in Y\times Z$. For background on absolute sums
the reader is referred to \cite{MPRY} and references therein.

\begin{proposition}[\mbox{\rm \cite[Proposici\'{o}n~1]{MarEMA}}] \label{prop:EMA}
Let $X$ be a Banach space and let $Y$, $Z$ be closed subspaces of
$X$. Suppose that $X$ is the absolute sum of $Y$ and $Z$. Then
$$
n(X)\leqslant \min\big\{n(Y),\ n(Z)\big\}.
$$
\end{proposition}

The following somehow surprising example was obtained in \cite{M-P}
by using Proposition~\ref{suma}.

\begin{example}[\mbox{\rm \cite[Example~2.b]{M-P}}] \label{exam-num0norma}
There is a real Banach space $X$ for which the numerical radius is a
norm but is not equivalent to the operator norm, i.e.\ the numerical
index of $X$ is $0$ although $v(T)>0$ for every non-null $T\in
L(X)$.
\end{example}

The numerical index of some vector-valued function spaces was also
computed in \cite{L-M-M,M-P,M-V}. Given a real or complex Banach
space $X$ and a compact Hausdorff topological space $K$, we write
$C(K,X)$ and $C_w(K,X)$ to denote, respectively, the space of
$X$-valued continuous (resp.\ weakly continuous) functions on $K$.
If $\mu$ is a positive $\sigma$-finite measure, by $L_1(\mu,X)$ and
$L_\infty(\mu,X)$ we denote respectively the space of $X$-valued
$\mu$-Bochner-integrable functions and the space of $X$-valued
$\mu$-Bochner-measurable functions which are essentially bounded.

\begin{theorem}[\mbox{\rm \cite{L-M-M},
\cite{M-P}, and \cite{M-V}}] Let $K$ be a compact Hausdorff space,
and let $\mu$ be a positive $\sigma$-finite measure. Then
$$
n(C_w(K,X))=n(C(K,X))=n(L_1(\mu,X))=n(L_\infty(\mu,X))=n(X)
$$
for every Banach space $X$.
\end{theorem}

The numerical index of $C_{w^*}(K,X^*)$, the space of $X^*$-valued
weakly-star continuous functions on $K$ is also studied in
\cite{L-M-M}. Unfortunately, only a partial result is achieved.

\begin{proposition}[\mbox{\rm \cite[Propositions 5 and 7]{L-M-M} }]
Let $K$ be a compact Hausdorff space and let $X$ be a Banach space.
Then
$$
n(C_{w^*}(K,X^*))\leqslant n(X).
$$
If, in addition, $X$ is an Asplund space or $K$ has a dense subset
of isolated points, then
$$
n(X^*)\leqslant n(C_{w^*}(K,X^*)).
$$
\end{proposition}

To finish this section let us comment that, roughly speaking, when
one finds an explicit computation of the numerical index of a Banach
space in the literature only few values appear; namely, $0$ (real
Hilbert spaces), $\e^{-1}$ (Glickfeld's example), $1/2$ (complex
Hilbert spaces), and $1$ ($C(K)$, $L_1(\mu)$,  and many more). The
preceding results about sums and vector-valued function spaces do
not help so much, and the exact values of $n(\ell_p^2)$ are not
still known. Let us also say that, when the authors of
\cite{D-Mc-P-W} prove \eqref{eq-DMCPW}, they only use examples of
Banach spaces whose numerical indices are the extremes of the
intervals, and then a connectedness argument is applied. Very
recently, M.~Mart\'{\i}n and J.~Mer\'{\i} have partially covered this gap in
\cite{MarMer}, where they explicitly compute the numerical index for
four families of norms on $\R^2$. The most interesting one is the
family of regular polygons.

\begin{proposition}[\mbox{\rm \cite[Theorem~5]{MarMer}}]\label{prop:regularpolygons}
Let $n\geqslant 2$ be a positive integer, and let $X_n$ be the
two-dimensional real normed space whose unit ball is the convex hull
of the $(2n)^{\rm th}$ roots of unity, i.e.\ $B_{X_n}$ is a regular
$2n$-polygon centered at the origin and such that one of its
vertices is $(1,0)$. Then,
$$
n(X_n)= \begin{cases}
 \tan\left(\dfrac{\pi}{2n}\right) & \emph{\text{if}} \ \ n \ \ \emph{\text{is even,}} \\
     &      \\
\sin\left(\dfrac{\pi}{2n}\right) & \emph{\text{if}} \ \ n \ \
\emph{\text{is odd.}}
\end{cases}
$$
\end{proposition}

\begin{problem}
Compute the numerical index for other families of finite-dimensional
Banach spaces. In particular, it would be interesting to get the
complex analog of the results in {\rm \cite{MarMer}}.
\end{problem}

\section{Numerical index and duality}\label{sec:duality}
Given a bounded linear operator $T$ on a Banach space $X$, it is a
well-known fact in the theory of numerical ranges (see \cite[\S
9]{B-D1}) that
\begin{equation}\label{eq:Lumer}
\sup \re V(T)=\lim_{\alpha\downarrow 0} \dfrac{\|\Id  + \alpha\, T\|
-1}{\alpha}
\end{equation}
and so
\begin{equation*}
v(T)=\max_{\omega\in \T}\,\lim_{\alpha\downarrow 0} \dfrac{\|\Id  +
\alpha\,\omega\, T\| -1}{\alpha}.
\end{equation*}
It is immediate from the above formula that
\begin{equation}\label{eq:numradius-adjoint}
v(T)=v(T^*),
\end{equation}
where $T^*\in L(X^*)$ is the adjoint operator of $T$, and the result
given in \cite[Proposition~1.3]{D-Mc-P-W} that
\begin{equation}\label{eq:inequality}
n(X^*)\leqslant n(X)
\end{equation}
clearly follows. The question if this is actually an equality had
been around from the beginning of the subject (see
\cite[pp.~386]{K-M-RP-survey}, for instance). Let us comment some
partial results which led to think that the answer could be
positive. Namely, it is clear that $n(X)=n(X^*)$ for every reflexive
space $X$, and this equality also holds whenever $n(X^*)=1$ , in
particular when $X$ is an $L$- or an $M$-space. Moreover, it is also
true that $n(X)=n(X^*)$ when $X$ is a $C^*$-algebra or a von Neumann
algebra predual (Theorem~\ref{th:cstaralgebras}).

Nevertheless, in a very recent paper \cite{BKMW}, K.~Boyko,
V.~Kadets, M.~Mart\'{\i}n, and D.~Werner have answered the question in
the negative by giving an example of a Banach space whose numerical
index is strictly greater than the numerical index of its dual. Our
aim in this section is to present such counterexample with a new and
more direct proof.

As usual, $c$ denotes the Banach space of all convergent scalar
sequences \mbox{$x=(x(1),x(2), \ldots)$} equipped with the sup-norm.
The dual space of $c$ is (isometric to) $\ell_1$ and we will write
$c^*\equiv \ell_1\oplus_1 \K$ where
$$
\left\langle (y,\lambda)\, ,\,x\right\rangle = \sum_{n=1}^\infty
y(n)\,x(n)\ +\ \lambda\,\lim x \qquad \big(x\in c,\ (y,\lambda)\in
\ell_1\oplus_1 \K\big).
$$
For every $n\in\mathbb{N}$, we denote by $e_n^*$ the norm-one
element of $c^*$ given by
$$
e_n^*(x)=x(n) \qquad \big(x\in c \big).
$$
We are now ready to show that the numerical index of a Banach space
and the one of its dual do not always coincide.

\begin{example}[\mbox{\rm \cite[Example~3.1]{BKMW}}] \label{example:main}
Let us consider the Banach space
$$
X=\big\{(x,y,z)\in c\oplus_\infty c \oplus_\infty c\ : \ \lim x +
\lim y + \lim z =0\big\}.
$$
Then, $n(X)=1$ and $n(X^*)<1$.
\end{example}

\begin{proof}
We observe that
\begin{align*}
X^*&=\bigl[c^*\oplus_1 c^* \oplus_1 c^*\bigr]/
\langle(\lim,\lim,\lim)\rangle
\end{align*}
so that, writing $Z=\ell_1^{3}/ \langle(1,1,1)\rangle$, we can
identify
\begin{equation}\label{eq:xstar-xdoublestar}
X^* \equiv \ell_1\oplus_1 \ell_1 \oplus_1 \ell_1 \oplus_1 Z \qquad
\text{and} \qquad  X^{**}  \equiv \ell_\infty \oplus_\infty
\ell_\infty \oplus_\infty \ell_\infty \oplus_\infty Z^*.
\end{equation}
\begin{figure}[h]
\centering \boxed{
\begin{minipage}[c]{80mm}
           \centering
        \resizebox{75mm}{!}{\includegraphics{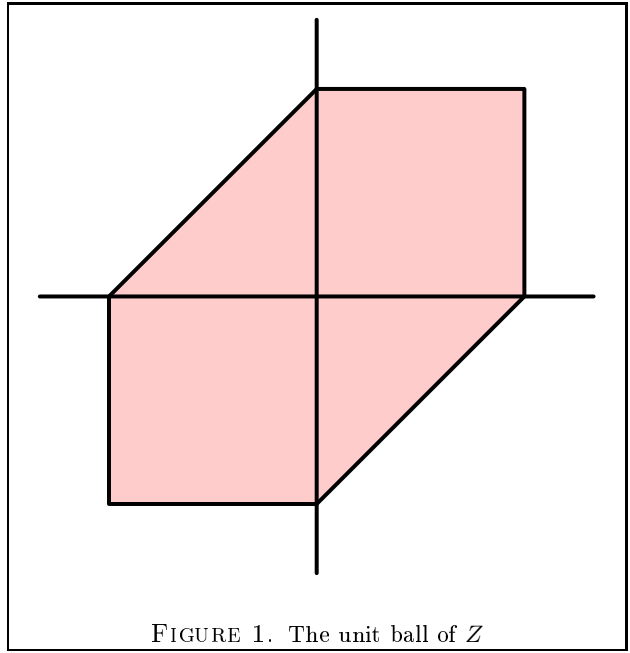}}
        \caption{\small The unit ball of $Z$}
        \label{figure}
            \end{minipage}}
\end{figure}

\noindent With this in mind, we write $A$ to denote the set
$$
\{(e_n^*,0,0,0)\, : \, n\in \mathbb{N}\}\cup \{(0,e_n^*,0,0)\, : \,
n\in \mathbb{N}\} \cup \{(0,0,e_n^*,0)\, : \, n\in \mathbb{N}\}.
$$
Then $A$ is clearly a norming subset of $S_{X^*}$ and
\begin{equation}\label{eq:extreme}
|x^{**}(a^*)|=1 \qquad \big(x^{**}\in \text{ext}(B_{X^{**}}), \
a^*\in A\big).
\end{equation}
Let us prove that $n(X)=1$. Indeed, we fix $T\in L(X)$ and $\eps>0$.
Since $T^*$ is $w^*$-continuous and $A$ is norming, we may find
$a^*\in A$ such that
$$
\|T^*(a^*)\|\geqslant \|T\|-\varepsilon.
$$
Now, we take $x^{**}\in \ext[X^{**}]$ such that
$$
|x^{**}(T^*(a^*))|=\|T^*(a^*)\|.
$$
Since $|x^{**}(a^*)|=1$ thanks to \eqref{eq:extreme}, we get
\begin{equation*}
v(T)=v(T^*)\geqslant |x^{**}(T^*(a^*))|\geqslant \|T\|-\varepsilon.
\end{equation*}
It clearly follows that $v(T)=\|T\|$ and $n(X)=1$.

To show that $n(X^*)<1$, we use \eqref{eq:xstar-xdoublestar} and
Proposition~\ref{suma} to get
$$
n(X^*)\leqslant n(Z),
$$
and the fact that $n(Z)<1$ follows easily from a result due to
C.~McGregor \cite[Theorem 3.1]{Mc}. Actually, in the real case, the
unit ball of $Z$ is an hexagon (see Figure~\ref{figure} above),
which is isometrically isomorphic to the space $X_3$ of
Proposition~\ref{prop:regularpolygons}, so $n(Z)=1/2$.\qed
\end{proof}

The above example can be pushed forward, to produce even more
striking counterexamples.

\begin{examples}[\mbox{\rm \cite[Examples~3.3]{BKMW}}]\ \label{examples:better} $ $
\begin{enumerate}
\item[(a)] There exists a real Banach space $X$
such that $n(X)=1$ and $n(X^*)=0$.
\item[(b)] There exists a complex Banach space $X$
satisfying that \mbox{$n(X)=1$} and $n(X^*)=1/\e$.
\end{enumerate}
\end{examples}

Once we know that the numerical index of a Banach space and the one
of its dual may be different, the question aries if two preduals of
a given Banach space have the same numerical index. By a
\emph{predual} of a Banach space $Y$ we mean a Banach space $X$ such
that $X^*$ is (isometrically isomorphic to) $Y$. The answer is again
negative, as the following result shows.

\begin{example}[\mbox{\rm \cite[Example~3.6]{BKMW}}]
\label{example:predual} Let us consider the Banach spaces
\begin{align*}
X_1&=\big\{(x,y,z)\in c\oplus_\infty c \oplus_\infty c\ : \ \lim x +
\lim y + \lim z =0\big\}\\
\intertext{and}
 X_2&=\big\{(x,y,z)\in c\oplus_\infty c \oplus_\infty
c\ : \ x(1) + y(1) + z(1) =0\big\}.
\end{align*}
Then, $X_1^*$ and $X^*_2$ are isometrically isomorphic, but
$n(X_1)=1$ and $n(X_2)<1$.
\end{example}

The following question might also be addressed.

\begin{problem}\label{problem:predualwiththesame}
Let $Y$ be a dual space. Does there exist a predual $X$ of $Y$ such
that $n(X)=n(Y)$?
\end{problem}

Another interesting issue could be to find isomorphic properties of
a Banach space $X$ ensuring that $n(X^*)=n(X)$. On the one hand,
Example~\ref{example:main} shows that Asplundness is not such a
property. On the other hand, it is shown in
\cite[Proposition~4.1]{BKMW} that if a Banach space $X$ with the
Radon-Nikod\'{y}m property has numerical index~$1$, then $X^*$ has
numerical index~$1$ as well. Therefore, the following question
naturally arises.

\begin{problem}
Let $X$ be a Banach space with the Radon-Nikod\'{y}m property. Is it
true that $n(X)=n(X^*)$?
\end{problem}

Another sufficient condition would follow from a positive answer to
Problem~\ref{problem:predualwiththesame}.

\begin{problem}
Let $Y$ be a dual space admitting a unique predual $X$ (up to
isometric isomorphisms). Is it true that $n(Y)=n(X)$?
\end{problem}

Let us finish this section by remarking that the space given in
Example~\ref{example:main} is useful as a counterexample for many
other conjectures, as we will see later on.

\section{Banach spaces with numerical index one}
The guiding open question on these spaces is the following.
\label{sec:numindex-1}

\begin{problem}
Find necessary and sufficient conditions for a Banach space to have
numerical index~$1$ which do not involve operators.
\end{problem}

In 1971, C.~McGregor \cite[Theorem 3.1]{Mc} gave such a
characterization in the finite-dimensional case. More concretely, a
finite-dimensional normed space $X$ has numerical index~$1$ if and
only if
\begin{equation}\label{eq:McGregor}
|x^*(x)|=1 \quad \text{for every $x\in \ext[X]$ and every $x^*\in
\ext[X^*]$.}
\end{equation}

It is not clear how to extend this result to arbitrary Banach
spaces. If we use literally \eqref{eq:McGregor} in the
infinite-dimensional context, we do not get a sufficient condition,
since the set $\ext[X]$ may be empty and this does not imply
numerical index~$1$ (e.g.\ $\ext[c_0(\ell_2)]=\emptyset$ but
$n(c_0(\ell_2))<1$). On the other hand, we do not know if
\eqref{eq:McGregor} is a necessary condition.

\begin{problem}
Let $X$ be a Banach space with numerical index $1$. Is it true that
$|x^*(x)|=1$ for every $x^*\in \ext[X^*]$ and every $x\in \ext[X]$?
\end{problem}

Our first aim in this section is to discuss several reformulations
of assertion \eqref{eq:McGregor} to get either sufficient or
necessary conditions for a Banach space to have numerical index~$1$.

Aiming at sufficient conditions, it is not difficult to show that
\eqref{eq:McGregor} implies numerical index~$1$ for a Banach space
$X$ as soon as the set $\ext[X]$ is large enough to determine the
norm of operators on $X$, i.e.\ $B_X=\ecc(\ext[X])$. Actually, we
may replace $\ext[X]$ with any subset of $S_X$ satisfying the same
property. On the other hand, we may replace $\ext[X]$ by
$\ext[X^{**}]$ and the role of $\ext[X^*]$ can be played by any
norming subset of $S_{X^*}$. Let us comment that this is is what we
did in the proof of Example~\ref{example:main}. All these ideas
appear implicity in several papers (see \cite{L-M-P,MarRNP,Mar-ADP}
for example); we summarize them in the following proposition.

\begin{proposition}\label{prop:normingimpliesnx1}
Let $X$ be a Banach space. Then, any of the following three
conditions is sufficient to ensure that $n(X)=1$.
\begin{enumerate}
\item[(a)] There exists a subset $C$ of $S_X$ such that
$\ecc(C)=B_X$ and
$$
|x^{*}(c)|=1
$$
for every $x^{*}\in\ext[X^{*}]$ and every $c\in C$.
\item[(b)] $|x^{**}(x^*)|=1\,$ for every $x^{**}\in \ext[X^{**}]$ and
every $x^*\in\ext[X^*]$.
\item[(c)] There exists a norming subset $A$ of $S_{X^*}$ such that
$$
|x^{**}(a^*)|=1
$$
for every $x^{**}\in\ext[X^{**}]$ and every $a^*\in A$.
\end{enumerate}
\end{proposition}

Let us comment on the converse of the above result. First, condition
(a) is not necessary as shown by $c_0$. Second, it was proved in
\cite[Example~3.4]{BKMW} that condition (b) is not necessary either,
the counterexample being the space given in
Example~\ref{example:main}. Finally, we do not know if there exists
a Banach with numerical index~$1$ in which condition (c) is not
satisfied.

\begin{problem}[\mbox{\rm \cite[Remark~3.5]{BKMW}}]
Let $X$ be a Banach space with numerical index~$1$. Does there exist
a norming subset $A$ of $S_{X^*}$ such that $|x^{**}(a^*)|=1$ for
every $x^{**}\in\ext[X^{**}]$ and every $a^*\in A$?
\end{problem}

Necessary conditions in the spirit of McGregor's result were given
in 1999 by G.~L\'{o}pez, M.~Mart\'{\i}n, and R.~Pay\'{a} \cite{L-M-P}. Their key
idea was considering denting points instead of general extreme
points. Recall that $x_0\in B_X$ is said to be a \emph{denting
point} of $B_X$ if it belongs to slices of $B_X$ with arbitrarily
small diameter. If $X$ is a dual space and the slices can be taken
to be defined by weak$^*$-continuous functionals, then we say that
$x_0$ is a \emph{weak$^*$-denting point}.

\begin{proposition}[\mbox{\rm \cite[Lemma~1]{L-M-P}}] \label{lemafundamental}
Let $X$ be a Banach space with numerical index~$1$. Then,
\begin{enumerate}
\item[(a)] $|x^*(x)|=1\ $ for every $x^*\in \ext[X^*]$ and every denting point $x\in
B_X$.
\item[(b)] $|x^{**}(x^*)|=1\ $ for every $x^{**}\in \ext[X^{**}]$ and every
weak$^*$-denting point \mbox{$x^*\in B_{X^*}$}.
\end{enumerate}
\end{proposition}

This result will play a key rolle in the next section.

Let us comment that, like McGregor original result, the conditions
in Proposition~\ref{lemafundamental} are not sufficient in the
infinite-dimensional context. Indeed, the space $X=C([0,1],\ell_2)$
does not have numerical index~$1$, while $B_X$ has no denting points
and there are no $w^*$-denting points in $B_{X^*}$. Actually, all
the slices of $B_X$ and the $w^*$-slices of $B_{X^*}$ have
diameter~$2$ (see \cite[Lemma~2.2 and Example on p.~858]{KSSW}, for
instance).

Anyhow, if we have a Banach space $X$ such that $B_X$ has enough
denting points (if $X$ has the Radon-Nikod\'{y}m property, for
instance), then item (a) in the above proposition combines with
Proposition~\ref{prop:normingimpliesnx1} to characterize the
numerical index~$1$ for $X$. The same is true for item (b) when
$B_{X^*}$ has enough weak$^*$-denting points (if $X$ is an Asplund
space, for instance).

\begin{corollary}[\mbox{\rm \cite[Theorem~1]{MarRNP} and \cite[\S 1]{Mar-ADP}}]
Let $X$ be a Banach space.
\begin{enumerate}
\item[(a)] If $X$ has the Radon-Nikod\'{y}m property, then the following
are equivalent:
\begin{enumerate}
\item[$(i)$] $X$ has numerical index~$1$.
\item[$(ii)$] $|x^*(x)|=1\ $ for every $x^*\in \ext[X^*]$ and every
denting point $x$ of $B_X$.
\item[$(iii)$] $|x^{**}(x^*)|=1\ $ for every $x^{**}\in \ext[X^{**}]$ and
every $x^*\in\ext[X^*]$.
\end{enumerate}
\vspace{1ex}
\item[(b)] If $X$ is an Asplund space, then the following are
equivalent:
\begin{enumerate}
\item[$(i)$] $X$ has numerical index~$1$.
\item[$(ii)$] $|x^{**}(x^*)|=1\ $ for every $x^{**}\in \ext[X^{**}]$ and every
weak$^*$-denting point \mbox{$x^*\in B_{X^*}$}.
\end{enumerate}
\end{enumerate}
\end{corollary}

In the remaining part of this section we discuss other kind of
sufficient conditions for a Banach space to have numerical
index~$1$.

The eldest of these properties was introduced in the fifties by
O.~Hanner \cite{Hann}: a real Banach space has the
\emph{intersection property 3.2} (\emph{3.2.I.P.} in short) if every
collection of three mutually intersecting closed balls has nonempty
intersection. The 3.2.I.P.\ was systematically studied by
J.~Lindenstrauss \cite{Lin} and \AA.~Lima \cite{Lima}, and typical
examples of spaces with this property are $L_1(\mu)$ and their
isometric preduals. Real Banach spaces with the 3.2.I.P.\ have
numerical index~$1$ since they fulfil
Proposition~\ref{prop:normingimpliesnx1}.b (see
\cite[Corollary~3.3]{Lima} and \cite[Theorem~4.7]{Lin}). The
converse is false even in the finite-dimensional case (see
\cite[Remark~3.6]{Hann} and \cite[p.~47]{Lin}).

Another isometric property, weaker than the 3.2.I.P.\ but still
ensuring numerical index~$1$, was introduced by R.~Fullerton in 1960
\cite{Full}. A real or complex Banach space is said to be a
\emph{CL-space} if its unit ball is the absolutely convex hull of
every maximal convex subset of the unit sphere. If the unit ball is
merely the closed absolutely convex hull of every maximal convex
subset of the unit sphere, we say that the space is an
\emph{almost-CL-space} (J.~Lindenstrauss \cite{Lin} and \AA.~Lima
\cite{Lima2}). Both definitions appeared only for real spaces, but
they extend literally to the complex case. For general information,
we refer to the already cited papers \cite{Lima,Lima2,Lin}; more
recent results can be found in \cite{MartPaya-CL,Reis}. Let us
remark that the complex space $\ell_1$ is an almost-CL-space which
is not a CL-space \cite[Proposition~1]{MartPaya-CL}, but we do not
know if such an example exists in the real case.

\begin{problem}
Is there any real almost-CL-space which is not a CL-space?
\end{problem}

The fact that CL-spaces have numerical index~$1$ was observed by
M.~Acosta \cite{Aco}, and her proof extends easily to
almost-CL-spaces (see \cite[Proposition~12]{Mar}). Actually,
almost-CL-spaces fulfil condition (c) of
Proposition~\ref{prop:normingimpliesnx1} as shown in
\cite[Lemma~3]{MartPaya-CL}. In the converse direction, the basic
examples of Banach spaces with numerical index~$1$ are known to be
almost-CL-spaces (see \cite{MartPaya-CL} and
\cite[Theorem~32.9]{B-D2}). Moreover, all finite-dimensional spaces
with numerical index~$1$ are CL-spaces \cite[Corollary~3.7]{Lima2},
and a Banach space with the Radon-Nikod\'{y}m property and numerical
index~$1$ is an almost-CL-space \cite[Theorem~1]{MarRNP}.
Nevertheless, Banach spaces with numerical index~$1$ which are no
almost-CL-spaces have been recently found. Actually, this happens
with the space given in Example~\ref{example:main}
\cite[Example~3.4]{BKMW}.

\begin{figure}[h]
\centering \boxed{
\begin{minipage}[c]{80mm}
           \centering
        \resizebox{75mm}{!}{\includegraphics{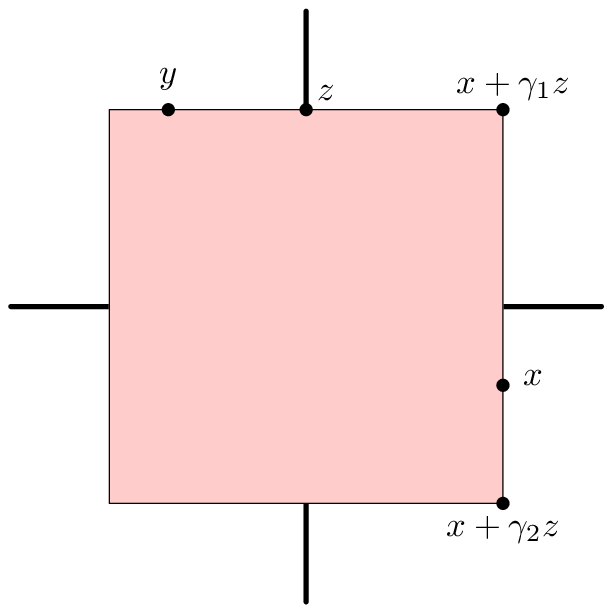}}
        \caption{\small $\ell_\infty^2$ is lush}
        \label{figure-lush}
            \end{minipage}}
\end{figure}

The last condition we would like to mention is a weakening of the
concept of almost-CL-space introduced in \cite{BKMW}. A Banach space
$X$ is said to be \emph{lush} if for every $x,y\in S_X$ and every
$\eps>0$, there exists $y^* \in S_{Y^*}$ such that
$$
y \in S(B_X,y^*, \eps):=\{z\in B_X\ : \ \re y^*(z)>1-\eps\}
$$
and
$$
\dist\bigl(x,\ec\bigl(\T\, S(B_X,y^*,\eps)\bigr)\bigr) < \eps.
$$
In the real case, the above definition is equivalent to the
following one: for every $x,y\in S_X$ and every $\eps>0$, there
exist $z\in S_X$ and $\gamma_1,\gamma_2\in\R$ such that
$$
\|y + z\|>2-\eps, \qquad |\gamma_1-\gamma_2|=2, \qquad
\|x+\gamma_i\,z\|\leqslant 1 + \eps \ \ (i=1,2).
$$
See Figure~\ref{figure-lush} below for an interpretation of this
property in dimension~$2$.

Almost-CL-spaces are clearly lush, and lush spaces have numerical
index~$1$ \cite[Proposition~2.2]{BKMW}. The converse of the first
implication is not true, and once again the counterexample is the
space $X$ given in Example~\ref{example:main}: as we already
mentioned, $X$ is not an almost-CL-space, but the original proof
given in \cite[Example~3.1]{BKMW} of the fact that $n(X)=1$ passes
through the lushness of the space. Actually, this example is only a
specimen of a general family of lush subspaces of $C(K)$ introduced
in \cite{BKMW}, namely C-rich subspaces. A subspace $X$ of $C(K)$ is
\emph{C-rich} if for every nonempty open subset $U$ of $K$ and every
$\eps > 0$, there is a positive continuous norm-one function $h$
with support inside $U$, such that the distance from $h$ to $X$ is
less than $\eps$.

\begin{theorem}[\mbox{\rm \cite[Theorem~2.4]{BKMW}}]
Let $K$ be a Hausdorff topological space and let $X$ be a C-rich
subspace of $C(K)$. Then, $X$ is lush and, therefore, $n(X)=1$.
\end{theorem}

In particular, this occurs for finite codimensional subspaces of
$C(K)$ when $K$ has no isolated points. Actually, the following
result characterizes C-rich finite-codimensional subspaces of an
arbitrary $C(K)$. We recall that the \emph{support} of a regular
measure $\mu\in C(K)^*$ is
$$
\text{supp}(\mu)=\bigcap\left\{C\subset K\ : \ C\ \text{closed},\
|\mu|(K\setminus C)=0\right\}.
$$

\begin{example}[\mbox{\rm \cite[Proposition~2.5]{BKMW}}]
Let $K$ be a compact Hausdorff space and $\mu_1,\ldots,\mu_n\in
C(K)^*$. The subspace
$$
Y=\bigcap_{i=1}^n  \ker \mu_i
$$
is C-rich if and only if $\ \bigcup_{i=1}^n \text{\rm supp}(\mu_i)$
does not intersect the set of isolated points of $K$.
\end{example}

We do not know if the class of lush spaces exhausts the whole class
of Banach spaces with numerical index~$1$.

\begin{problem}
Let $X$ be a Banach space with $n(X)=1$. Is $X$ lush?
\end{problem}

\section{Renorming and numerical index}\label{sec:renorming}
In 2003, C.~Finet, M.~Mart\'{\i}n, and R.~Pay\'{a} \cite{F-M-P} studied the
numerical index from the isomorphic point of view, i.e.\ they
investigated the set $\mathcal{N}(X)$ of those values of the
numerical index which can be obtained by equivalent renormings of a
Banach space $X$. This study has a precedent in the 1974 paper
\cite{Tille} by K.~Tillekeratne, where it is proved that every
complex space of dimension greater than one can be renormed to
achieve the minimum value of the numerical index; the same is true
for real spaces.

\begin{proposition}[\mbox{\rm \cite[Proposition~1]{F-M-P} and
\cite[Theorem~3.1]{Tille}}] \label{prop:cero1/e} Let $X$ be a Banach
space of dimension greater than one. Then $\,0\in \mathcal{N}(X)$ in
the real case, $\e^{-1}\in \mathcal{N}(X)$ in the complex case.
\end{proposition}

One of the main aims of \cite{F-M-P} is to show that
$\mathcal{N}(X)$ is a interval for every Banach space $X$. To get
this result, the authors use the continuity of the mapping carrying
every equivalent norm on $X$ to its numerical index with respect to
a metric taken from \cite[\S 18]{B-D2}.

\begin{proposition}[\mbox{\rm \cite[Proposition~2]{F-M-P}}]\label{prop:continuity}
$\mathcal{N}(X)$ is an interval for every Banach space $X$.
\end{proposition}

As an immediate consequence of the above two results, we get the
following.

\begin{corollary}[\mbox{\rm \cite[Corollary~3]{F-M-P}}]\label{cor:1implieseverything}
If $\,1\in \mathcal{N}(X)$ for a Banach space $X$ of dimension
greater than one, then $\mathcal{N}(X)=[0,1]$ in the real case and
$\mathcal{N}(X)=[\e^{-1},1]$ in the complex case.
\end{corollary}

Since $n(\ell_\infty^m)=1$ for every $m$, the following particular
case arises.

\begin{corollary}[\mbox{\rm \cite[Theorem 3.2]{Tille}}]
Let $m$ be an integer larger than $1$. Then
$$
\mathcal{N}(\R^m)=[0,1] \qquad \text{ and } \qquad
\mathcal{N}(\C^m)=[\e^{-1},1].
$$
\end{corollary}

Now, one may ask if the above result is also true in the
infinite-dimensional context, equivalently, whether or not every
Banach space can be equivalently renormed to have numerical
index~$1$. The answer is negative, as shown in the already cited
paper \cite{L-M-P}.

\begin{theorem}[\mbox{\rm \cite[Theorem~3]{L-M-P}}]\label{th:LMP-containsl1}
Let $X$ be an infinite-dimensional real Banach space with $1\in
\mathcal{N}(X)$. If $X$ has the Radon-Nikod\'{y}m property, then $X$
contains $\ell_1$. If $X$ is an Asplund space, then $X^*$ contains
$\ell_1$.
\end{theorem}

It follows that infinite-dimensional real reflexive spaces cannot be
renormed to have numerical index~$1$. But even more is true.

\begin{corollary}[\mbox{\rm \cite[Corollary~5]{L-M-P}}]\label{cor:nonseparable}
Let $X$ be an infinite-dimensional real Banach space. If  $X^{**}/X$
is separable, then $1\notin \mathcal{N}(X)$.
\end{corollary}

It is easy to explain how Theorem~\ref{th:LMP-containsl1} was proved
in \cite{L-M-P}. Namely, the authors used
Proposition~\ref{lemafundamental}, the well-known facts that the
unit ball of a space with the Radon-Nikod\'{y}m property has many
denting points and that the dual unit ball of an Asplund space has
many weak$^*$-denting points (see \cite{Bou}, for instance), and the
following sufficient condition for a real Banach space to contain
either $c_0$ or $\ell_1$.

\begin{lemma}[\mbox{\rm \cite[Proposition~2]{L-M-P}}] Let $X$ be a
real Banach space, and assume that there is an infinite set
$A\subset S_X$ such that $|x^*(a)| = 1$ for every $a\in A$ and every
$x^*\in \ext[X^*]$. Then $X$ contains $c_0$ or $\ell_1$.
\end{lemma}

Thus, the first open question in this line is the following.

\begin{problem}
Characterize those Banach spaces which can be equivalently renormed
to have numerical index $1$.
\end{problem}

We propose to study separately necessary and sufficient conditions
for a Banach space to be renormable with numerical index~$1$. With
respect to necessary conditions, we have obtained two in the real
case, namely Theorem~\ref{th:LMP-containsl1} and
Corollary~\ref{cor:nonseparable}. It is not known if they are valid
in the complex case; actually, the following especial case remains
open.

\begin{problem}
Does there exist an infinite-dimensional complex reflexive space
which can be renormed to have numerical index~$1$?
\end{problem}

For more necessary conditions, we suggest to study the following two
questions. The first one has a positive answer for real
almost-CL-spaces \cite[Theorem~5]{MartPaya-CL}.

\begin{problem}\label{prob:duall_1}
Let $X$ be an infinite-dimensional (real) Banach space satisfying
that $1\in \mathcal{N}(X)$. Does $X^*$ contain $\ell_1$?
\end{problem}

\begin{problem}
Let $X$ be an infinite-dimensional (real) Banach space satisfying
that $1\in \mathcal{N}(X)$. Does $X$ contain $c_0$ or $\ell_1$?
\end{problem}

Let us remark that we do not know of any non-trivial sufficient
condition for a Banach space to be renormable with numerical
index~$1$. We propose the following ones to be checked.

\begin{problem}
Let $X$ be a Banach space containing an infinite-dimensional
subspace $Y$ such that $1\in\mathcal{N}(Y)$. Is it true that
$1\in\mathcal{N}(X)$?
\end{problem}

One may consider many especial cases.

\begin{problem}
Let $X$ be a Banach space containing a subspace isomorphic to either
$c_0$, $\ell_1$, $C[0,1]$, or $L_1[0,1]$. Is it true that
$1\in\mathcal{N}(X)$?
\end{problem}

The following question is especially interesting since in view of
Problem~\ref{prob:duall_1} it might lead to a characterization of
Banach spaces that can be renormed to have numerical index~$1$.

\begin{problem}
Let $X$ be a Banach space such that $\ell_1\subseteq X^*$. Is is
true that $1\in\mathcal{N}(X)$?
\end{problem}

We finish this section by showing that the value $1$ of the
numerical index is very particular. Indeed, it is proved in
\cite{F-M-P} that ``most'' Banach spaces can be renormed to achieve
any possible value for the numerical index except eventually $1$.
Recall that a system $\{(x_\lambda,x^*_\lambda)\}_{\lambda\in
\Lambda}\subset X\times X^*$ is said to be \emph{biorthogonal} if
$x_\lambda^*(x_\mu)=\delta_{\lambda,\mu}$ for $\lambda,\mu\in
\Lambda$, and \emph{long} if the cardinality of $\Lambda$ coincides
with the density character of $X$.

\begin{theorem}[\mbox{\rm \cite[Theorem~10]{F-M-P}}]\label{th:alpha}
Let $X$ be a Banach space admitting a long biorthogonal system. Then
$\sup\mathcal{N}(X)=1$. Therefore, when the dimension of $X$ is
greater than one, $\mathcal{N}(X)\supset [0,1[$ in the real case and
$\mathcal{N}(X)\supset [e^{-1},1[$ in the complex case.
\end{theorem}

Typical examples of Banach spaces admitting a long biorthogonal
system are WCG spaces (see \cite{DGZ}). For instance, if $X^{**}/X$
is separable, then the Banach space $X$ is WCG (see
\cite[Theorem~3]{Val}, for example) while, in the real case,
$1\notin \mathcal{N}(X)$ unless $X$ is finite-dimensional (see
Corollary~\ref{cor:nonseparable}). Therefore, in many cases one of
the inclusions of Theorem~\ref{th:alpha} becomes an equality.

\begin{corollary}[\mbox{\rm \cite[Corollary~11]{F-M-P}}]
Let $X$ be an infinite-dimensional real Banach space such that
$X^{**}/X$ is separable. Then $\mathcal{N}(X)=[0,1[$.
\end{corollary}

Let us comment that Theorem~\ref{th:alpha} is proved by using a
geometrical property that was introduced by J.~Lindenstrauss in the
study of norn-attaining operators \cite{Lin0} and called
\emph{property $\alpha$} by W.~Schachermayer \cite{Scha}. It is
known that, under the continuum hypothesis, there are Banach spaces
which cannot be renormed with property $\alpha$ \cite{G-T, Neg}.
Nevertheless, B.~Godun and S.~Troyanski proved in
\cite[Theorem~1]{G-T} that this renorming is possible for Banach
spaces admitting a long biorthogonal system; as far as we know, this
is the largest class of spaces for which renorming with property
$\alpha$ is possible.

The question arises if the assumption of having a long biorthogonal
system in Theorem~\ref{th:alpha} can be dropped.

\begin{problem}
Is is true that $\sup \mathcal{N}(X)=1$ for every Banach space $X$?
\end{problem}

It is also studied in \cite{F-M-P} the relationship between the
numerical index and the so-called \emph{property $\beta$}
\cite{Lin0,Scha}. Contrary to property $\alpha$, property $\beta$ is
isomorphically trivial (J.~Partington \cite{Par}), but it does not
produce such a good result as Theorem~\ref{th:alpha}. At least, it
can be used to prove that $\mathcal{N}(X)$ does not reduces to a
point when the dimension of $X$ is greater than one.

\begin{theorem}[\mbox{\rm \cite[Theorem~9]{F-M-P}}]\label{th:beta}
Let $X$ be a Banach space with $\dim(X)>1$. Then
$\mathcal{N}(X)\supset [0,1/3[$ in the real case and
$\mathcal{N}(X)\supset [\e^{-1},1/2[$ in the complex case.
\end{theorem}

\section{Real Banach spaces with numerical index zero}
\label{sec:numindex0}  As we commented in the introduction, real
Banach spaces underlying complex Banach spaces as well as real
Hilbert spaces of dimension greater than one, have numerical
index~$0$. By Proposition~\ref{prop:EMA}, the absolute sum of such a
space and any other real Banach space has numerical index~$0$.

The general open question in this section is the following.

\begin{problem}\label{prob:ni0-char}
Find characterizations of Banach spaces with numerical index~$0$
which do not involve operators.
\end{problem}

A sufficient condition which generalizes all the introductory
examples is the following easy result of M.~Mart\'{\i}n, J.~Mer\'{\i} and
A.~Rodr\'{\i}guez Palacios \cite{MaMeRo}. We say that a real vector space
has a \emph{complex structure} if it is the real space underlying a
complex vector space.

\begin{proposition}[\mbox{\rm \cite[Proposition~2.1]{MaMeRo}}]\label{prop-suff}
Let $X$ be a real Banach space, and let $Y, Z$ be closed subspaces
of $X$, with $Z\neq 0$. Suppose that $Z$ is endowed with a complex
structure, that $X=Y \oplus Z$, and that the equality $\left\|y +
\e^{i\rho} z \right\| =\|y+ z\|$ holds for every $(\rho ,y,z)\in \R
\times Y\times Z$. Then $n(X)=0$.
\end{proposition}

It looks like if a necessary condition for numerical index~$0$ could
be the emergence of a subspace with some kind of complex structure.
As a matter of fact, the following example shows that this is not
the case.

\begin{example}[\mbox{\rm \cite[Example~2.2]{MaMeRo}}]\label{example:numindexzero-poly}
There exists a real Banach space $X$ with numerical index~$0$ which
is polyhedral, i.e.\ the intersection of $B_X$ with any
finite-dimensional subspace of $X$ is the convex hull of a finite
set of points. Therefore, $X$ does not contain any isometric copy of
$\C$.
\end{example}

The above example has the additional interest that the numerical
radius is a norm on $L(X)$, i.e.\ the only operator with numerical
radius $0$ is the zero operator. It could be the case that a Banach
space in which there is a non-null operator with numerical radius
$0$ has a subspace with some kind of complex structure. We recall
that a bounded linear operator $T$ on a (real or complex) Banach
space $X$ is \emph{skew-hermitian} if $\re V(T)= \{0\}$; we write
$\mathcal{Z}(X)$ for the (possibly null) closed subspace of $L(X)$
consisting of skew-hermitian operators on $X$. In the real case, $T$
is skew-hermitian if and only if $v(T)=0$; when the space $X$ is
complex, an operator $T$ is \emph{hermitian} if $V(T)\subset \R$,
i.e.\ the operator $i\,T$ is skew-hermitian. Hermitian operators
have been deeply studied since the sixties and many results on
Banach algebras depend on them; we refer to \cite{B-D1,B-D2} for
more information. Also, skew-hermitian operators have been widely
discussed in the seventies and eighties, especially in the
finite-dimensional case; more information can be found in the papers
by H.~Rosenthal \cite{Rosenthal,Ros-Pacific} and references therein.

\begin{problem}\label{prob:ni0-complex-structure}
Let $X$ be a real space which has a non-null skew-hermitian
operator. Does $X$ contain a subspace with a complex structure?
\end{problem}

Let us give a clarifying example. If $H$ is a $n$-dimensional
Hilbert space, it is easy to check that $\mathcal{Z}(H)$ is the
space of skew-symmetric operators on $H$ (i.e.\ $T^*=-T$ in the
Hilbert space sense), so it identifies with the space of
skew-symmetric matrixes $\mathcal{A}(n)$. It is a classical result
from the theory of linear algebra that a $n\times n$ matrix $A$
belongs to $\mathcal{A}(n)$ if and only if $\exp(\rho A)$ is an
orthogonal matrix for every $\rho\in\R$ (see
\cite[Corollary~8.5.10]{Artin-Algebra} for instance).

It is shown in \cite[\S3]{B-D1} that the above fact extends to
general Banach spaces. Indeed, for an arbitrary Banach space $X$ and
an operator $T\in L(X)$, by making use of the ``exponential
formula''
$$
\sup \re V(T)=\sup_{\alpha>0}\ \frac{\log\|\exp(\alpha
T)\|}{\alpha}\, ,
$$
it is easy to prove that the following are equivalent:
\begin{enumerate}
\item[$(i)$] $T$ is skew-hermitian,
\item[$(ii)$] $\exp(\rho T)$ is an onto isometry for every $\rho\in \R$.
\end{enumerate}

From now on, we will restrict ourselves to finite-dimensional
spaces. Here, the above result reminds the theory of Lie groups and
Lie algebras, for which the group of orthogonal matrices and the
space of skew-symmetric matrices are distinguished examples.
Actually, the group of all the isometries on a finite-dimensional
real space $X$ is a Lie group whose associated Lie algebra is
$\mathcal{Z}(X)$ \cite[Theorem~1.4 and Proposition~1.5]{Rosenthal}.
With this in mind and using results from the theory of Lie groups,
it is proved in \cite[Theorem~3.8]{Rosenthal} that the following are
equivalent for a finite-dimensional real space $X$:
\begin{enumerate}
\item[$(i)$] the numerical index of $X$ is $0$,
\item[$(ii)$] there are infinitely many isometries on $X$.
\end{enumerate}

The main open question in this section is the following.

\begin{problem}\label{prob:describe-fin-dim-ni0}
Describe the finite-dimensional real Banach spaces with numerical
index~$0$.
\end{problem}

The next result follows this line. In particular, it shows that
finite-dimensional normed spaces with numerical index~$0$ wear some
kind of complex structure.

\begin{theorem}[\mbox{\rm \cite[Corollary~3.7]{Rosenthal} and \cite[Theorem~2.4]{MaMeRo}}]
\label{th:numindexzero} Let $X$ be a finite-dimensional real Banach
space. Then, the following are equivalent:
\begin{enumerate}
 \item[$(i)$] The numerical index of $X$ is zero.
 \item[$(ii)$] There are nonzero complex vector spaces
$X_1,\ldots,X_m$, a real vector space $X_0$, and positive integer
numbers $q_1,\dots,q_m$ such that $X=X_0\oplus X_1\oplus \dots
\oplus X_m$ and
$$
\left\|x_0+\e^{iq_1\rho}x_1+\dots+\e^{iq_m\rho}x_m\right\|=
\|x_0+x_1+\dots+x_m\|
$$
for all $\rho \in \R$, $x_j \in X_j$ $(j=0,1,\dots,m)$.
\end{enumerate}
\end{theorem}

Some remarks are pertinent. First, the above result shows that a
finite-dimensional space with numerical index~$0$ contains complex
subspaces (at least one), whose complex structures are well related
one to the others and to the non-complex part. Second, the above
result is taken literally from the paper \cite{MaMeRo}, but the same
equivalence with arbitrary real numbers $q_1,\ldots,q_m$ had
appeared in the 1985 Rosenthal's paper \cite{Rosenthal}; the fact
that the $q_j$'s can be taken integers is new from \cite{MaMeRo} and
it uses the classical Kronecker's Approximation Theorem (see
\cite[Theorem~442]{H-W} for instance). Third, it can be asked if the
number of complex spaces in the theorem can be reduced to one, and
so Proposition~\ref{prop-suff} would be an equivalence in the
finite-dimensional case. The answer is negative.

\begin{example}[\mbox{\rm \cite[Example~2.8]{MaMeRo}}]
The space $\R^4$ with norm
$$
\|(a,b,c,d)\|=\dfrac14 \int_0^{2\pi} \left|\re \left(\e^{2i t}(a + i
b) + \e^{it}(c + id)\right)\right|\, dt \qquad (a,b,c,d\in\R)
$$
has numerical index~$0$, but the number of complex spaces in
Theorem~\ref{th:numindexzero}.$ii$ cannot be reduced to one.
\end{example}

Such an example as the above is not possible in dimensions two or
three. Actually, in this case, Theorem~\ref{th:numindexzero} takes a
more suitable form \cite[Theorem~3.1]{Rosenthal}. Let $X$ be a real
Banach space with numerical index~$0$.
\begin{enumerate}
\item[(a)] If $\dim(X)=2$, then $X$ is isometrically isomorphic to
the two-dimensional real Hilbert space.
\item[(b)] If $\dim(X)=3$, then $X$ is an absolute sum of $\R$ and the
two-dimensional real Hilbert space.
\end{enumerate}

Our next aim is to discuss some questions related to the Lie algebra
of skew-hermitian operators $\mathcal{Z}(X)$ of an arbitrary
$n$-dimensional space which. The main related open question is the
following.

\begin{problem}\label{prob:dimzx}
Figure out what are the possible values for the dimension of
$\mathcal{Z}(X)$ when $\dim(X)=n$.
\end{problem}

Let us fix a $n$-dimensional Banach space $X$. It follows from a
theorem of Auerbach \cite[Theorem~9.5.1]{Rolewicz}, that there
exists an inner product $(\cdot|\cdot)$ on $X$ such that every
skew-hermitian operator on $X$ remains skew-hermitian (hence
skew-symmetric) on $H:=(X,(\cdot|\cdot))$. Then, by just fixing an
orthonormal basis of $H$, we get an identification of
$\mathcal{Z}(X)$ with a Lie subalgebra of the Lie algebra
$\mathcal{A}(n)$. Therefore,
$$
\dim\big(\mathcal{Z}(X)\big)\leqslant \dfrac{n(n-1)}{2}.
$$
The equality holds if and only if $X$ is a Hilbert space (see
\cite[Theorem~3.2]{Rosenthal} or \cite[Corollary~2.7]{MaMeRo}). It
is a good question whether or not all the intermediate numbers are
possible values for the dimension of $\mathcal{Z}(X)$. The answer is
negative, as a consequence of Theorem~3.2 in Rosenthal's paper
\cite{Rosenthal}, which reads as follows.
\begin{enumerate}
\item[(a)] If $\dim \big(Z(X)\big)>\frac{(n-1)(n-2)}{2}$, then $X$ is a
Hilbert space and so $\dim \big(Z(X)\big)=\frac{n(n-1)}{2}$.
\item[(b)] $\dim\big(Z(X)\big)=\frac{(n-1)(n-2)}{2}$ if and only if
$X$ is a non-Euclidean absolute sum of $\R$ and a Hilbert space of
dimension $n-1$.
\end{enumerate}
For low dimensions, Problem~\ref{prob:dimzx} has been solved in
\cite{Ros-Pacific}. When the dimension of $X$ is $3$, the above
result leaves only the following possible values for the dimension
of $\mathcal{Z}(X)$: $0$ as for $X=\ell_\infty^3$, $1$ as for
$\R\oplus_1 \C$, and $3$ as for $\ell_2^3$. When the dimension of
$X$ is $4$, the possible values of the dimension of $\mathcal{Z}(X)$
allowed by the above result are $0,1,2,3,6$; all of them are
possible \cite[pp.~443]{Ros-Pacific}. The first dimension in which
Problem~\ref{prob:dimzx} is open is $n=5$.

\begin{problem}
What are the possible values for the dimension of $\mathcal{Z}(X)$
when $X$ is a $5$-dimensional real Banach space?
\end{problem}

\section{Asymptotic behavior of the set of finite-dimensional
spaces with numerical index one}\label{sec:fin-dim-numindex1} Let us
start the section by recalling that most of the sufficient and the
necessary conditions for having numerical index~$1$ given in
section~$3$ are actually characterizations in the finite-dimensional
case. For the convenience of the reader, we summarize them in the
following paragraph.
\begin{center}
\begin{minipage}[c]{120mm}\label{page:finite-dimension}
{\slshape Let $X$ be a finite-dimensional real normed space. Then,
the following are equivalent.
\begin{enumerate}
\item[$(i)$] $X$ has numerical index~$1$.
\item[$(ii)$] $|x^*(x)|=1$ for all the extreme points $x^*$ of
$B_{X^*}$ and $x$ of $B_X$.
\item[$(iii)$] $X$ is a CL-space, i.e.~$B_X=\ec(F\cup -F)$ for every
maximal convex subset $F$ of $S_X$.
\item[$(iv)$] $X$ is lush, i.e.\ for every $x,y\in S_X$, there are
$z\in S_X$ and $\gamma_1,\gamma_2\in \R$ such that $\|y + z\|=2$,
$\,|\gamma_1-\gamma_2|=2$, and $\|x+\gamma_i\,z\|\leqslant 1 \ \
(i=1,2)$.
\end{enumerate}
 }
\end{minipage}
\end{center}

Our aim in this section is to consider the asymptotic behavior (as
the dimension grows to infinity) of some parameters related to the
Banach-Mazur distance for the family of finite-dimensional real
normed spaces with numerical index~$1$. Let us write $\mathcal{N}_m$
for the space of all $m$-dimensional normed spaces endowed with the
\emph{Banach-Mazur distance}
$$
d(X,Y)=\inf\{\|T\|\,\|T^{-1}\|\ : \ T:X\longrightarrow Y \text{
isomorphism}\} \qquad (X,Y\in \mathcal{N}_m),
$$
and let us write $\mathcal{M}_m$ for the subset consisting of those
$m$-dimensional spaces with numerical index~$1$. Our aim is to study
some questions related to these two spaces. As far as we know, the
first result of this kind was given very recently by T.~Oikhberg
\cite{Oik04}.

\begin{theorem}[\mbox{\rm \cite[Theorem~4.1]{Oik04}}] There exists a
universal positive constant $c$ such that
$$
d(X,\ell_2^m)\geqslant c\,m^{\frac14}
$$
for every $m\geqslant 1$ and every $X\in \mathcal{M}_m$.
\end{theorem}

It is well-known that
$d(\ell_1^m,\ell_2^m)=d(\ell_\infty^m,\ell_2^m)=\sqrt{m}$ for every
$m>1$ (see \cite[pp.~720]{Gia-Milman} for instance). Therefore, the
following question arises naturally.

\begin{problem}
Does there exists a universal constant $c>0$ such that
$$
d(X,\ell_2^m)\geqslant c\,\sqrt{m}
$$
for every $m\geqslant 1$ and every $X\in\mathcal{M}_m$?
\end{problem}

It was observed in \cite[pp.~622]{Oik04} that the answer to this
question is positive for the spaces with the 3.2.I.P., even with
$c=1$. This is a consequence of the 1981 result by A.~Hansen and
\AA.~Lima \cite{H-L} that these spaces are constructed starting from
the real line and producing successively $\ell_\infty$ and/or
$\ell_1$ sums. But, as we already mentioned, not every element of
$\mathcal{M}_m$ has the 3.2.I.P.

Finally, we would like to propose some related questions.

\begin{problem}
What is the diameter of $\mathcal{M}_m$? Is it (asymptotically)
close to the diameter of $\mathcal{N}_m$?
\end{problem}

\begin{problem}
What is the biggest possible distance from an element of
$\mathcal{N}_m$ to the set $\mathcal{M}_m$?
\end{problem}

\section{Relationship to the Daugavet
property.}\label{sec:ADP}

In every Banach space with the Radon-Nikod\'{y}m property (in particular
in every reflexive space) the unit ball must have denting points.
There are Banach spaces $X$ (as $C[0,1]$, $L_1[0,1]$, and many
others) with an extremely opposite property: for every $x\in S_X$
and for arbitrarily small $\eps> 0$, the closure of
$$
\ec\bigl(B_X \setminus (x + (2 -\eps)B_X)\bigr)
$$
equals to the whole $B_X$ (see Figure~\ref{figure-Daugavet} below).
\begin{figure}[h,t] \centering \boxed{
\begin{minipage}[c]{80mm}
           \centering
        \resizebox{75mm}{!}{\includegraphics{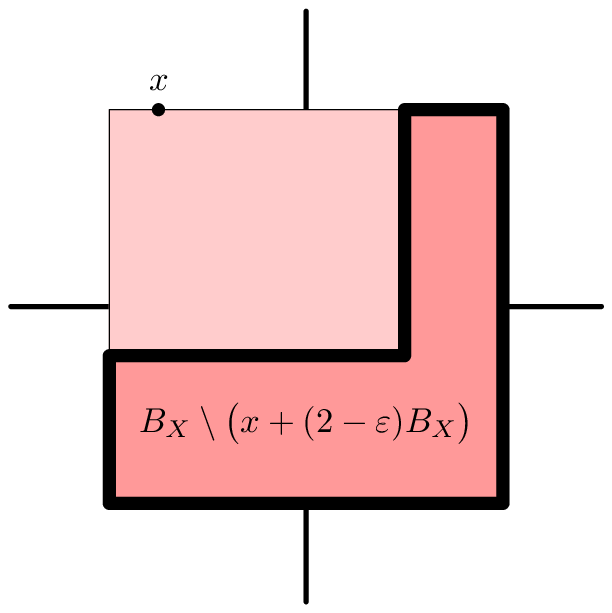}}
        \caption{\small The space $\ell_\infty^2$ does not have the Daugavet property}
        \label{figure-Daugavet}
            \end{minipage}}
\end{figure}
This geometric property of the space is equivalent to the following
exotic property of operators on $X$: for every compact operator
$T:X\longrightarrow X$, the so-called \emph{Daugavet equation}
\begin{equation}\label{DE}\tag{\textrm{DE}}
\|\Id + T \| = 1 + \|T\|
\end{equation}
holds. This property of $C[0,1]$ was discovered by I.~Daugavet in
1963 and is called the \emph{Daugavet property} \cite{KSSW0,KSSW}.
Over the years, the validity of the Daugavet equation was proved for
some classes of operators on various spaces, including weakly
compact operators on $C(K)$ and $L_1(\mu)$ provided that $K$ is
perfect and $\mu$ does not have any atoms (see \cite{Wer0} for an
elementary approach), and on certain function algebras such as the
disk algebra $A(\mathbb{D})$ or the algebra of bounded analytic
functions $H^\infty$ \cite{WerJFA,Woj}. In the nineties, new ideas
were infused into this field and the geometry of Banach spaces
having the Daugavet property was studied; we cite the papers of
V.~Kadets, R.~Shvidkoy, G.~Sirotkin, and D.~Werner \cite{KSSW} and
R.~Shvidkoy \cite{Shv} as representatives. Let us comment that the
original definition of Daugavet property given in \cite{KSSW0,KSSW}
only required rank-one operators to satisfy \eqref{DE} and, in such
a case, this equation also holds for every bounded operator which
does not fix a copy of $\ell_1$ \cite{Shv}.

Although the Daugavet property is of isometric nature, it induces
various isomorphic restrictions. For instance, a Banach space with
the Daugavet property contains $\ell_1$ \cite{KSSW}, it does not
have unconditional basis (V.~Kadets \cite{Kadets}) and, moreover, it
does not isomorphically embed into an unconditional sum of Banach
spaces without a copy of $\ell_1$ \cite{Shv}. It is worthwhile to
remark that the latest result continues a line of generalization
(\cite{K97}, \cite{KS99}, \cite{KSSW}) of the well known theorem by
A.~Pe\l czy\'nski \cite{Pelcz} that $L_1[0,1]$ (and so $C[0,1]$)
does not embed into a space with unconditional basis.

The state-of-the-art on the Daugavet property can be found in
\cite{WerSur}; for very recent results we refer the reader to
\cite{BM,BKSW,I-K-W,KaKaWe,KaWe} and references therein.

Let us explain the relation between \eqref{DE} and the numerical
range of an operator. In the aforementioned paper \cite{D-Mc-P-W} by
J.~Duncan, C.~McGregor, J.~Pryce, and A.~White, it was deduced from
formula \eqref{eq:Lumer} on page~\pageref{eq:Lumer} that an operator
$T$ on a Banach space $X$ satisfies \eqref{DE} if and only if $\sup
\re V(T)=\|T\|$. Therefore, $v(T)=\|T\|$ if and only the following
equality holds
\begin{equation}\label{aDE}\tag{\textrm{aDE}}
\max_{\omega\in \T} \|\Id  + \omega\,T\|= 1 + \|T\|
\end{equation}
(see \cite[Lemma~2.3]{MaOi} for an explicit proof). Therefore, it
was known since 1970 that every bounded linear operator on $C(K)$ or
$L_1(\mu)$ satisfies \eqref{aDE}, a fact that was rediscovered and
reproved in some papers from the eighties and nineties as the ones
by Y.~Abramovich \cite{Abra0}, J.~Holub \cite{Holub}, and K.~Schmidt
\cite{Sch}.

This latest equation was named as the \emph{alternative Daugavet
equation} by M.\ Mart\'{\i}n and T.\ Oikhberg in \cite{MaOi}, where the
following property was introduced. A Banach space $X$ is said to
have the \emph{alternative Daugavet property} if every rank-one
operator on $X$ satisfies \eqref{aDE}. In such a case, every weakly
compact operator on $X$ also satisfies \eqref{aDE}
\cite[Theorem~2.2]{MaOi}. Therefore, $X$ has the alternative
Daugavet property if and only if $v(T)=\|T\|$ for every weakly
compact operator $T\in L(X)$.

Let us comment that, contrary to the Daugavet property, this
property depends upon the base field (e.g.\ $\C$ has it as a complex
space but not as a real space). For more information on the
alternative Daugavet property we refer to the already cited paper
\cite{MaOi} and also to \cite{Mar-ADP}. From the former one we take
the following geometric characterizations of the alternative
Daugavet property.

\begin{proposition}[\mbox{\rm \cite[Propositions 2.1 and 2.6]{MaOi}}]
\label{prop:ADPcharac} Let $X$ be a Banach space. Then, the
following are equivalent.
\begin{enumerate}
\item[$(i)$] $X$ has the alternative Daugavet property.
\item[$(ii)$] For all $x_0\in S_X$, $x_0^*\in S_{X^*}$ and
$\eps>0$, there is some $x\in S_X$ such that
$$
|x_0^*(x)|\geqslant 1 - \eps \qquad \text{and} \qquad \|x +
x_0\|\geqslant 2-\eps.
$$
\item[$(ii^*)$] For all $x_0\in S_X$, $x_0^*\in S_{X^*}$ and
$\eps>0$, there is some $x^*\in S_{X^*}$ such that
$$
|x^*(x_0)|\geqslant 1 - \eps \qquad \text{and} \qquad \|x^* +
x_0^*\|\geqslant 2-\eps.
$$
\item[$(iii)$] $\displaystyle B_{X}=\ecc\Bigl(\T\,\bigl[B_X \setminus \bigl(x + (2
-\eps)B_X\bigr)\bigr]\Bigr)$ for every $x\in S_X$ and every $\eps>0$
(see Figure~\ref{figure-ADP} below).
\item[$(iii^*)$] $\displaystyle B_{X^*}=\ecc^{w^*}\Bigl(\T\,\bigl[B_{X^*} \setminus \bigl(x^* + (2
-\eps)B_{X^*}\bigr)\bigr]\Bigr)$ for every $x^*\in S_{X^*}$ and
every $\eps>0$.
\item[(iv)] $\displaystyle B_{X^*\oplus_\infty X^{**}} = \ecc^{w^*}\big(\{(x^*,x^{**})\, : \ x^*\in\ext[X^*],\,
x^{**}\in \ext[X^{**}],\ |x^{**}(x^*)|=1\}\big)$.
\end{enumerate}
\end{proposition}
\begin{figure}[h,t] \centering \boxed{
\begin{minipage}[c]{80mm}
           \centering
        \resizebox{75mm}{!}{\includegraphics{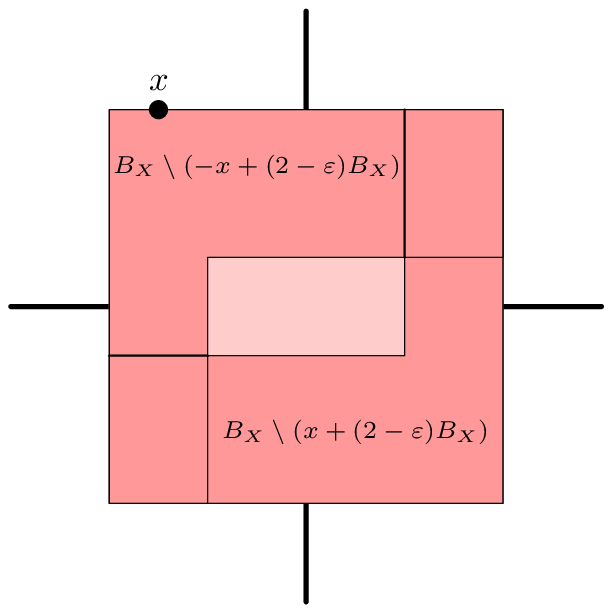}}
        \caption{\small The space $\ell_\infty^2$ has the alternative Daugavet property}
        \label{figure-ADP}
            \end{minipage}}
\end{figure}

It is clear that both spaces with the Daugavet property and spaces
with numerical index~$1$ have the alternative Daugavet property.
Both converses are false: the space $c_0 \oplus_1 C([0,1],\ell_2)$
has the alternative Daugavet property but fails the Daugavet
property and its numerical index is not $1$
\cite[Example~3.2]{MaOi}. Nevertheless, under certain isomorphic
conditions, the alternative Daugavet property forces the numerical
index to be $1$.

\begin{proposition}[\mbox{\rm \cite[Remark~6]{L-M-P}}] \label{prop:RNP-Asplund-ADP=>n(X)=1}
Let $X$ be a Banach space with the alternative Daugavet property. If
$X$ has the Radon-Nikod\'{y}m property or $X$ is an Asplund space, then
$n(X)=1$.
\end{proposition}

With this result in mind, one realizes that the necessary conditions
for a real Banach space to be renormed with numerical index~$1$
given in section~$4$ (namely Theorem~\ref{th:LMP-containsl1} and
Corollary~\ref{cor:nonseparable}), can be written in terms of the
alternative Daugavet property. Even more, in the proof of
Proposition~\ref{lemafundamental} given in \cite{L-M-P}, only
rank-one operators are used and, therefore, it can be also written
in terms of the alternative Daugavet property.

\begin{proposition}[\mbox{\rm \cite[Lemma~1 and Remark~6]{L-M-P}}]
\label{lemma:lemafundamentalADP} Let $X$ be a Banach space with the
alternative Daugavet property. Then,
\begin{enumerate}
\item[(a)] $|x^{**}(x^*)|=1\ $ for every $x^{**}\in \ext[X^{**}]$ and every
weak$^*$-denting point $x^*\in B_{X^*}$.
\item[(b)] $|x^*(x)|=1\ $ for every $x^*\in \ext[X^*]$ and every denting point $x\in
B_X$.
\end{enumerate}
\end{proposition}

\begin{proposition}[\mbox{\rm \cite[Remark~2.8]{MaOi}}]
Let $X$ be an infinite-dimensional real Banach space with the
alternative Daugavet property. If $X$ has the Radon-Nikod\'{y}m
property, then $X$ contains $\ell_1$. If $X$ is an Asplund space,
then $X^*$ contains $\ell_1$. In particular, $X^{**}/X$ is not
separable.
\end{proposition}

The above two results give us an indication of why it is difficult
to find characterizations of Banach spaces with numerical index~$1$
that do no involve operators. Indeed, it is not easy to construct
noncompact operators on an abstract Banach space. Thus, when on uses
the assumption that a Banach space has numerical index~$1$, only the
alternative Daugavet property can be easily exploited. Of course,
things are easier if one is working in a context where the
alternative Daugavet property ensures numerical index~$1$, as it
happens with Asplund spaces and spaces with the Radon-Nikod\'{y}m
property. Therefore, it would be desirable to find more isomorphic
properties ensuring that the alternative Daugavet property implies
numerical index~$1$. One possibility is the following.

\begin{problem}
Let $X$ be a Banach space which does not contain any copy of
$\ell_1$ and having the alternative Daugavet property. Is it true
that $n(X)=1$?
\end{problem}

An affirmative answer to the above question might come from the
following more ambitious one.

\begin{problem}
Let $X$ be a Banach space with the alternative Daugavet property. Is
it true that $v(T)=\|T\|$ for every operator $T\in L(X)$ fixing no
copy of $\ell_1$?
\end{problem}

On the other hand, it is not possible to find isomorphic properties
ensuring that the alternative Daugavet property and the Daugavet
property are equivalent.

\begin{proposition}[\mbox{\rm \cite[Corollary~3.3]{MaOi}}]
Let $X$ be a Banach space with the alternative Daugavet property.
Then there exists a Banach space $Y$, isomorphic to $X$, which has
the alternative Daugavet property but fails the Daugavet property.
\end{proposition}

We may then look for isometric conditions that allow passing from
the alternative Daugavet property to the Daugavet property. Having a
complex structure could be such a condition.

\begin{problem}
Let $X$ be a complex Banach space such that $X_\R$ has the
alternative Daugavet property. Does it follow that $X$ (equivalently
$X_\R$) has the Daugavet property?
\end{problem}

\section{The polynomial numerical indices}\label{sec:poly}
In 1968, the concept of numerical range of operators was extended to
arbitrary continuous functions from the unit sphere of a real or
complex Banach space into the space by F.~Bonsall, B.~Cain, and
H.~Schneider \cite{BonCaiSch} in the obvious way. They showed that
the numerical range of a bounded linear operator is always
connected, and the same is true for continuous functions with only
one exception when the space is real and of dimension one. Three
years later, L.~Harris \cite{Harris} studied various possible
numerical ranges for holomorphic functions on complex Banach spaces
and got some deep results on them.

In this section, we concentrate on homogeneous polynomials on real
or complex spaces. Let us give the necessary definitions; we refer
the reader to the book \cite{Dineenbook} by S.~Dineen for background
on polynomials. Given a Banach space $X$ and a positive integer $k$,
a mapping $P:X\longrightarrow X$ is called a (continuous)
\emph{\mbox{$k$-homogeneous} polynomial} on $X$ if there is an
$k$-linear continuous mapping $A:X\times \cdots \times X
\longrightarrow X$ such that $P(x)=A(x,\ldots,x)$ for every $x\in
X$. A \emph{polynomial} on $X$ is a finite sum of homogeneous
polynomials. The space $\mathcal{P}(X;X)$ of all polynomials on $X$
is normed by
$$
\|P\|=\sup_{x\in B_X}\,\|P(x)\| \qquad \bigl(P\in \mathcal{P}(X;
X)\bigr).
$$
We write $\mathcal{P}(^k X;X)$ to denote the subspace of
$\mathcal{P}(X;X)$ of those $k$-homogeneous polynomials on $X$,
which is a Banach space. The numerical range of $P\in
\mathcal{P}(X;X)$ is the set of scalars
$$
V(P):=\{x^*(P(x))\ : \ x\in S_X,\ x^*\in S_{X^*},\ x^*(x)=1\},
$$
and the numerical radius of $P$ is
$$
v(P):=\sup\{|\lambda|\ : \ \lambda\in V(P)\}.
$$
As in the linear case, it is natural to define the \emph{polynomial
numerical index of order $k$} of $X$ to be the constant
\begin{align*}
n^{(k)}(X)& :=\inf\left\{v(P)\ : \ P\in\mathcal{P}(^k X;X),
\|P\|=1\right\} \\
& = \max\left\{M\geqslant 0 \ : \ \|P\|\leqslant M v(P)\ \forall
P\in \mathcal{P}(^k X;X)\right\};
\end{align*}
of course, $n^{(1)}(X)$ coincides with the usual numerical index of
the space $X$. This definition was introduced very recently by
Y.~Choi, D.~Garc\'{\i}a, S.~Kim, and M.~Maestre \cite{ChoiGarcKimMaes}.
Let us remark that $0\leqslant n^{(k)}(X)\leqslant 1$, and that
$n^{(k)}(X)>0$ if and only if $v(\,\cdot\,)$ is a norm on
$\mathcal{P}(^k X;X)$ equivalent to the usual norm.

The first result that we would like to mention here is an extension
of Glickfeld's result for linear operators \cite{Gli} given by
L.~Harris \cite[Theorem~1]{Harris}: for complex spaces, the
numerical radius is always an equivalent norm in the space of
$k$-homogeneous polynomials. More concretely, if $X$ is a complex
Banach space and $k\geqslant 2$, then
\begin{equation}\label{eq:Harris}
n^{(k)}(X)\geqslant \exp\left(\frac{k\, \log(k)}{1-k}\right).
\end{equation}
It was also proved in \cite[\S7]{Harris} that the above inequalities
are sharp. In the real case, Harris' result above is false since the
polynomial numerical indices of a real space may vanish (see
Example~\ref{exa:num-index-k}.b below).

The next result shows the behavior of the polynomial numerical index
when the degree grows.

\begin{proposition}[\mbox{\rm
\cite[Proposition~2.5]{ChoiGarcKimMaes}}]\label{prop:numindex-k-decreasing}
Let $X$ be a Banach space and let $k$ be a positive integer. Then
$$
n^{(k+1)}(X)\leqslant n^{(k)}(X).
$$
\end{proposition}

Next we list several examples of Banach spaces for which there is
some information on their polynomial numerical indices.

\begin{examples}[\mbox{\rm \cite{ChoiGarcKimMaes}, \cite{ChoiGarcMartMaes}
and \cite{ChKi}}]\label{exa:num-index-k} $ $
\begin{enumerate}
\item[(a)] $n^{(k)}(\K)=1$ for every $k\in\N$.
\item[(b)] If $H$ is a real Hilbert space of dimension greater than
$1$, then $n^{(k)}(H)=0$ for every $k\in\N$.
\item[(c)] If $H$ is a complex Hilbert space of dimension greater
than $1$, then $1/4\leqslant n^{(k)}(H)\leqslant 1/2$ for every
$k\in\N$.
\item[(d)] In the complex case, all the polynomial numerical indices of
an $L_1$-predual are equal to $1$. In particular, this is the case
for the complex spaces $c_0$, $\ell_\infty$ and, more generally,
$C(K)$.
\item[(e)] $n^{(2)}(\ell_1)\leqslant 1/2$ in the real as well as in the complex
case.
\item[(f)] For every $k\geqslant 2$, the real spaces $c_0$, $\ell_\infty$, $c$, and
$\ell_\infty^m$ with $m\geqslant 2$, have numerical index of order
$k$ smaller than $1$.
\item[(g)] Actually, if $K$ is a non-perfect compact Hausdorff space with at
least two points, then the real space $C(K)$ has numerical index of
order $k$ smaller than $1$ for every $k\geqslant 2$.
\end{enumerate}
\end{examples}

The above examples show that the polynomial numerical indices
distinguish between $L$-spaces and $M$-spaces in the complex case,
and this is not possible if we only use the usual (linear) numerical
index. They also show that, unlike the linear case, there is no
relationship between the polynomial numerical indices of a Banach
space and the ones of its dual (see items (d) and (e) above).
Nevertheless, there is a relationship with the polynomial numerical
indices of the bidual.

\begin{proposition}[\mbox{\rm \cite[Corollary~2.15]{ChoiGarcKimMaes}}]
Let $X$ be a Banach space. Then, we have
$$
n^{(k)}(X^{**})\leqslant n^{(k)}(X)
$$
for every $k\geqslant 2$.
\end{proposition}

We do not know if the above inequalities are actually equalities.

\begin{problem}
Is there a Banach space $X$ such that $n^{(2)}(X^{**})< n^{(2)}(X)$?
\end{problem}

Let us present now more open questions arising from
Examples~\ref{exa:num-index-k}.

\begin{problem}
Compute the polynomial numerical indices of complex Hilbert spaces
and of $\ell_1$.
\end{problem}

\begin{problem}
Is $n^{(2)}\big(C[0,1]\big)=1$ in the real case? Is
$n^{(2)}\big(L_1[0,1]\big)=1$ in the real and/or in the complex
case?
\end{problem}

The behavior of polynomial numerical indices under direct sums and
the indices of vector-valued continuous functions spaces were also
studied in \cite{ChoiGarcKimMaes}. The results can be summarizedas
follows.

\begin{proposition}[\mbox{\rm \cite[Propositions 2.8 and 2.10]{ChoiGarcKimMaes}}]
\label{prop-numindexk-sum-vector} Let $k$ be a positive integer.
\begin{enumerate}
\item[(a)] If $\{X_\lambda:\ \lambda \in \Lambda\}$ is any family of Banach
spaces, then
\begin{align*}
n^{(k)}\Bigl(\left[\oplus_{\lambda\in\Lambda} X_\lambda\right]_{c_0}
\Bigr)&\leqslant \inf_{\lambda} \,n^{(k)}(X_\lambda) \\
n^{(k)}\Bigl(\left[\oplus_{\lambda\in\Lambda}
X_\lambda\right]_{\ell_1} \Bigr) &\leqslant \inf_{\lambda}
\,n^{(k)}(X_\lambda) \\
n^{(k)}\Bigl(\left[\oplus_{\lambda\in\Lambda}
X_\lambda\right]_{\ell_\infty} \Bigr)&\leqslant \inf_{\lambda}
\,n^{(k)}(X_\lambda).
\end{align*}
\item[(b)] If $X$ is a Banach space and $K$ is a compact Hausdorff space, then
$$
n^{(k)}\bigl(C(K,X)\bigr)\leqslant n^{(k)}(X).
$$
\end{enumerate}
\end{proposition}

Examples~\ref{exa:num-index-k} show that the inequalities in item
(a) above are not always equalities in the real case; just take
$\Lambda=\N$ and $X_n=\R$ for every $n\in\N$. In the complex case,
the same is true for the inequality involving $\ell_1$-sums
($\Lambda=\N$ and $X_n=\C$ for every $n\in\N$), but we do not know
the answer for the other sums.

\begin{problem}
Let $\{X_\lambda:\ \lambda \in \Lambda\}$ be a family of complex
Banach spaces. Is is true that
\begin{equation*}
n^{(k)}\Bigl(\left[\oplus_{\lambda\in\Lambda} X_\lambda\right]_{c_0}
\Bigr)= n^{(k)}\Bigl(\left[\oplus_{\lambda\in\Lambda}
X_\lambda\right]_{\ell_\infty} \Bigr) = \inf_{\lambda}
\,n^{(k)}(X_\lambda)\ ?
\end{equation*}
\end{problem}

For item (b) of Proposition~\ref{prop-numindexk-sum-vector} the
situation is similar. For instance, the real space $C(\{1,2\},\R)$
does not have the same polynomial numerical indices as $\R$. We do
not know if there exists an example of the same kind in the complex
case.

\begin{problem}
Let $X$ be a complex Banach space. Is $n^{(k)}\bigl(C(K,X)\bigr)$
equal to $n^{(k)}(X)$?
\end{problem}

For spaces of Bochner-measurable integrable or essentially bounded
vector valued functions no results are yet available.

\begin{problem}
Let $X$ be a Banach space, $\mu$ a positive $\sigma$-finite measure
and $k\geqslant 2$. Is there any relationship between
$n^{(k)}\bigl(L_1(\mu,X)\bigr)$ (resp.\
$n^{(k)}\bigl(L_\infty(\mu,X)\bigl)$) and $n^{(k)}(X)$?
\end{problem}

Item (a) in Proposition~\ref{prop-numindexk-sum-vector} can be used
to produce the following nice example.

\begin{example}
There exists a complex Banach space $X$ such that
$$
n^{(k)}(X)= \exp\left(\frac{k\, \log(k)}{1-k}\right)
$$
for every $k\geqslant 2$, i.e.\ all the inequalities in
\eqref{eq:Harris} are simultaneously equalities.
\end{example}

\begin{proof}
For each positive integer $k\geqslant 2$, let $X_k$ be the
two-dimensional complex Banach space given in \cite[\S7]{Harris}
which satisfies the required equality for this $k$. Then, the space
$$
X=\left[\oplus_{k\geqslant 2} X_k\right]_{c_0}
$$
satisfies
$$
n^{(k)}(X)\leqslant n^{(k)}(X_k)=\exp\left(\frac{k\,
\log(k)}{1-k}\right)
$$
by Proposition~\ref{prop-numindexk-sum-vector}.a, and the other
inequality always holds.\qed
\end{proof}

After presenting the known results and open questions related to the
computing of the polynomial numerical indices, we would like to
discuss the isometric or structural consequences that these indices
may have on a Banach space. Actually, we do not know of any result
in this line, so we simply state conjectures and open questions.

First, we do not know what are precisely the values that the
polynomial numerical index may take.

\begin{problem}
For a given $k\geqslant 2$, describe the sets
$$
\Bigl\{n^{(k)}(X)\ : \ X \text{ real Banach space}\Bigr\} \qquad
\text{and} \qquad \Bigl\{n^{(k)}(X)\ : \ X \text{ complex Banach
space}\Bigr\}.
$$
\end{problem}

We may also ask for the class of Banach spaces for which the the
$k$-order numerical index is one of the extreme values.

\begin{problem}
Study the real Banach spaces $X$ satisfying $n^{(2)}(X)=0$. For
instance, do they satisfy that $n(X)=0$? What is the answer for
finite-dimensional spaces?
\end{problem}

\begin{problem}
Characterize the complex Banach spaces $X$ satisfying $n^{(k)}(X)=1$
for all $k\geqslant 2$. In particular, do they always contain a copy
of $c_0$ in the infinite-dimensional case?
\end{problem}

\begin{problem}
Is there any real Banach space $X$ different from $\R$ such that
$n^{(2)}(X)=1$?
\end{problem}

It would be desirable to get more information on the non-increasing
sequence $\{n^{(k)}(X)\}_{k\in\N}$ of all polynomial numerical
indices of a Banach space $X$. For instance, in view of the
Examples~\ref{exa:num-index-k}, the following question arises.

\begin{problem}
Is there any Banach space $X$ for which $\displaystyle
\lim_{k\to\infty}\,n^{(k)}(X)\neq 0,1$?
\end{problem}

To finish the section, let us mention that Y.~Choi, D.~Garc\'{\i}a,
M.~Mart\'{\i}n, and M.~Maestre have extended the study of the Daugavet
equation to polynomials in the very recent paper
\cite{ChoiGarcMartMaes}. Given a Banach space $X$, a polynomial
$P\in \mathcal{P}(X;X)$ satisfies the Daugavet equation if
$$
\|\Id + P\|=1 + \|P\|,
$$
and it satisfies the alternative Daugavet equation if
$$
\max_{\omega\in \T} \|\Id  + \omega\,P\|= 1 + \|P\|.
$$
Daugavet and alternative Daugavet properties are also translated in
\cite{ChoiGarcMartMaes} to the polynomial setting by direct
generalization of the linear case. A Banach space $X$ has the
\emph{$k$-order Daugavet property} (resp.\ the \emph{$k$-order
alternative Daugavet property}) if every rank-one $k$-homogeneous
polynomial on $X$ satisfies the Daugavet equation (resp.\ the
alternative Daugavet equation). These properties are related to the
polynomial numerical indices since, as in the linear case, for a
polynomial $P$ we have that
\begin{enumerate}
\item[(a)] $P$ satisfies the Daugavet equation if and only if $\|P\|=\sup\re V(P)$,
\item[(b)] $P$ satisfies the alternative Daugavet equation if and only if
$\|P\|=v(P)$
\end{enumerate}
(see \cite[Proposition~1.3]{ChoiGarcMartMaes}). Let us remark that
part of the information given in Examples~\ref{exa:num-index-k}
comes from results on the Daugavet and the alternative Daugavet
properties \cite{ChoiGarcMartMaes}.

We conclude this paper by mentioning the following rather surprising
result.

\begin{proposition}[\mbox{\rm \cite[Proposition~3.3 and Remark~3.6]{ChoiGarcMartMaes}}]
Let $X$ be a Banach space and let $k\geqslant 2$. Then, the Daugavet
equation and the alternative Daugavet equations are equivalent in
$\mathcal{P}(^k X;X)$ in the following two cases:
\begin{enumerate}
\item[(a)] if the base field is $\C$, or
\item[(b)] if the base field is $\R$ and $k$ is even.
\end{enumerate}
If $k$ is odd, then the Daugavet and the alternative Daugavet
equation are not equivalent in $\mathcal{P}(^k \R;\R)$.
\end{proposition}

\bigskip

\bigskip

\noindent\textbf{Acknowledgments:\ } The work of the first named
author was partially supported by a fellowship from the
\emph{Alexander-von-Humboldt Stiftung}. Second and third authors are
partially supported by Spanish MCYT project no.\ BFM2003-01681 and
Junta de Andaluc\'{\i}a grant FQM-185.


\begin{thebibliography}{99}\small

\bibitem{Aco} \textsc{M.~D.~Acosta},
Operators that attain its numerical radius and CL-spaces,
\emph{Extracta Math.} \textbf{5} (1990), 138--140.

\bibitem{Abra0} \textsc{Y.~Abramovich},
A generalization of a theorem of J.~Holub, \emph{Proc. Amer. Math.
Soc.} \textbf{108} (1990), 937--939.

\bibitem{Artin-Algebra} \textsc{M.~Artin}, \emph{Algebra},
Prentice-Hall, Englewood Cliffs, New Jersey, 1991.

\bibitem{Bauer} \textsc{F.~L.~Bauer},
On the field of values subordinate to a norm, \emph{Numer. Math.}
\textbf{4} (1962), 103--111.

\bibitem{BM} \textsc{J.~Becerra Guerrero and M.~Mart\'{\i}n}, The
Daugavet property of $C^*$-algebras, $JB^*$-triples, and of their
isometric preduals, \emph{J. Funct. Anal.} \textbf{224} (2005),
316--337.

\bibitem{BKSW} \textsc{D.~Bilik, V.~Kadets, R.~V.~Shvidkoy, and
D.~Werner}, Narrow operators and the Daugavet property for
ultraproduts, \emph{Positivity} \textbf{9} (2005), 45--62.

\bibitem{B-K} \textsc{H.~F.~Bohnenblust and S.~Karlin},
Geometrical properties of the unit sphere in Banach algebras,
\emph{Ann. of Math.} \textbf{62} (1955), 217--229.

\bibitem{BonCaiSch} \textsc{F.~F.~Bonsall, B.~E.~Cain, and
H.~Schneider}, The numerical range of a continuos mapping of a
normed space, \emph{Aequationes Math.} \textbf{2} (1968), 86--93.

\bibitem{B-D1} \textsc{F.~F.~Bonsall and J.~Duncan},
\emph{Numerical Ranges of Operators on Normed Spaces and of Elements
of Normed Algebras}, London Math. Soc. Lecture Note Series
\textbf{2}, Cambridge, 1971.

\bibitem{B-D2} \textsc{F.~F.~Bonsall and J.~Duncan},
\emph{Numerical Ranges II}, London Math. Soc. Lecture Note Series
\textbf{10}, Cambridge, 1973.

\bibitem{Bou} \textsc{R.~R.~Bourgin},
\emph{Geometric Aspects of Convex Sets with the Radon-Nikod\'{y}m
Property}, Lecture Notes in Math. \textbf{993}, Springer-Verlag,
Berlin 1983.

\bibitem{BKMW} \textsc{K.~Boyko, V.~Kadets, M.~Mart\'{\i}n, and D.~Werner},
Numerical index of Banach spaces and duality, \emph{Math. Proc.
Cambridge Philos. Soc.} (to appear).

\bibitem{ChoiGarcKimMaes} \textsc{Y.~S.~Choi, D.~Garc\'{\i}a, S.~G.~Kim, and
M.~Maestre}, The polynomial numerical index of a Banach space,
\emph{Proc. Edinburgh Math. Soc.} \textbf{49} (2006), 39--52.

\bibitem{ChoiGarcMartMaes} \textsc{Y.~S.~Choi, D.~Garc\'{\i}a, M.~Mart\'{\i}n, and
M.~Maestre}, The Daugavet equation for polynomials, \emph{preprint}.

\bibitem{ChKi} \textsc{Y.~S.~Choi and S.~G.~Kim},
Norm or numerical radius attaining multilinear mappings and
polynomials, \emph{J. London Math. Soc.} \textbf{54} (1996),
135--147.

\bibitem{C-D-Mc} \textsc{M.~J.~Crabb, J.~Duncan, and C.~M.~McGregor},
Mapping theorems and the numerical radius, \emph{Proc. London Math.
Soc.} \textbf{25} (1972), 486--502.

\bibitem{Dau} \textsc{I.~K.~Daugavet},
On a property of completely continuous operators in the space $C$,
\emph{Uspekhi Mat. Nauk} \textbf{18} (1963), 157--158 (Russian).


\bibitem{DGZ} \textsc{R.~Deville, G.~Godefroy, and V.~Zizler},
\emph{Smoothness and Renormings in Banach spaces}, Pitman Monographs
64, New York 1993.

\bibitem{Dineenbook} \textsc{S. Dineen}, \emph{Complex Analysis on Infinite
Dimensional  Spaces}, Springer-Verlag, Springer Monographs in
Mathematics, Springer-Verlag, London,  1999.

\bibitem{D-Mc-P-W} \textsc{J.~Duncan, C.~McGregor, J.~Pryce, and
A.~White}, The numerical index of a normed space, \emph{J. London
Math. Soc.} \textbf{2} (1970), 481--488.

\bibitem{Ed-Dari} \textsc{E.~Ed-Dari},
On the numerical index of Banach spaces, \emph{Linear Algebra Appl.}
403 (2005), 86--96.

\bibitem{Ed-Dari-Khamsi} \textsc{E.~Ed-Dari and M.~A.~Khamsi},
The numerical index of the $L_p$ space, \emph{Proc. Amer. Math.
Soc.} \textbf{134} (2006), 2019--2025.

\bibitem{F-M-P} \textsc{C.~Finet, M.~Mart\'{\i}n, and R.~Pay\'{a}},
Numerical index and renorming, \emph{Proc. Amer. Math. Soc.}
\textbf{131} (2003), 871--877.

\bibitem{Full} \textsc{R.~E.~Fullerton},
Geometrical characterization of certain function spaces. \emph{In:
Proc. Inter. Sympos. Linear spaces (Jerusalem 1960)}, pp. 227--236.
Pergamon, Oxford 1961.

\bibitem{Gia-Milman} \textsc{A.~A.~Giannopoulos and V.~D.~Milman},
Euclidean structure in finite dimensional normed spaces,
\emph{Handbook of the geometry of Banach spaces, Vol. I},
pp.~707--779, North-Holland, Amsterdam, 2001.

\bibitem{Gli} \textsc{B.~W.~Glickfeld},
On an inequality of Banach algebra geometry and semi-inner-product
space theory, \emph{Illinois J. Math.} \textbf{14} (1970), 76--81.

\bibitem{G-T} \textsc{B.~V.~Godun and S.~L.~Troyanski},
Renorming Banach spaces with fundamental biorthogonal system,
\emph{Contemporary Math.} \textbf{144} (1993), 119--126.

\bibitem{GusRao} \textsc{K.~E.~Gustafson, and D.~K.~M.~Rao},
\emph{Numerical range. The field of values of linear operators and
matrices}, Springer-Verlag, New York, 1997.

\bibitem{Hal} \textsc{P.~Halmos},
\emph{A Hilbert space problem book}, Van Nostrand, New York, 1967.

\bibitem{Hann} \textsc{O.~Hanner},
Intersections of translates of convex bodies, \emph{Math. Scan.}
\textbf{4} (1956), 65--87.

\bibitem{H-L} \textsc{A.~B.~Hansen and \AA.~Lima},
The structure of finite dimensional Banach spaces with the 3.2.\
intersection property, \emph{Acta Math.} \textbf{146} (1981), 1--23.

\bibitem{H-W} \textsc{G.~H.~Hardy and E.~M.~Wright},
\emph{An introduction to the theory of numbers (5th edition)},
Clarendon Press, Oxford, 1979.

\bibitem{Harris} \textsc{L.~A.~Harris},
The numerical range of holomorphic functions in Banach spaces,
\emph{American J. Math.} \textbf{93} (1971), 1005--1019.

\bibitem{Holub} \textsc{J.~R.~Holub},
A property of weakly compact operators on $C[0,1]$, \emph{Proc.
Amer. Math. Soc.} \textbf{97} (1986), 396--398.

\bibitem{Hur} \textsc{T.~Huruya}, The normed space numerical index
of $C^*$-algebras, \emph{Proc. Amer. Math. Soc.} \textbf{63} (1977),
289--290.

\bibitem{IsiRod} \textsc{J.~M.~Isidro and A.~Rodr\'{\i}guez},
On the definition of real $W^*$-algebras, \emph{Proc. Amer. Math.
Soc.} \textbf{124} (1996), 3407--3410.

\bibitem{I-K-W} \textsc{Y.~Ivaknho, V.~Kadets, and D.~Werner}, The
Daugavet property for spaces of Lipschitz functions, \emph{Math.
Scand.} (to appear).

\bibitem{Kadets} \textsc{V.~M.~Kadets},
Some remarks concerning the Daugavet equation, \emph{Quaestiones
Math.} \textbf{19} (1996), 225--235.

\bibitem{K97} \textsc{V.~Kadets}, A generalization of a Daugavet's
theorem with applications to the space $C$ geometry,
\emph{Funktsional. Analiz i ego Prilozhen.} \textbf{31} (1997),
74--76. (Russian)

\bibitem{KaKaWe} \textsc{V.~Kadets, N.~Kalton, and D.~Werner}, Remarks on
rich subspaces of Banach spaces, \emph{Studia Math.} \textbf{159}
(2003), 195--206.

\bibitem{KS99} \textsc{V.~Kadets and R.~Shvidkoy}, The Daugavet
property for pairs of Banach  spaces, \emph{Matematicheskaya Fizika,
Analiz, Geometria}, \textbf{6} (1999), 253--263.

\bibitem{KSSW0} \textsc{V.~M.~Kadets, R.~V.~Shvidkoy, G.~G.~Sirotkin, and
D.~Werner}, Espaces de Banach ayant la propri\'{e}t\'{e} de Daugavet,
\emph{C. R. Acad. Sci. Paris}, S\'{e}r. I, \textbf{325} (1997),
1291--1294.

\bibitem{KSSW} \textsc{V.~Kadets, R.~Shvidkoy, G.~Sirotkin, and
D.~Werner}, Banach spaces with the Daugavet property, \emph{Trans.
Amer. Math. Soc.} \textbf{352} (2000), 855--873.

\bibitem{KaWe} \textsc{V.~Kadets and D.~Werner}, A Banach space
with the Schur and the Daugavet property, \emph{Proc. Amer. Math.
Soc.} \textbf{132} (2004), 1765--1773.

\bibitem{K-M-RP} \textsc{A.~Kaidi, A.~Morales, and
A.~Rodr\'{\i}guez-Palacios}, Geometrical properties of the product of a
$C^*$-algebra, \emph{Rocky Mountain J. Math.} \textbf{31} (2001),
197--213.

\bibitem{K-M-RP-survey} \textsc{A.~Kaidi, A.~Morales, and
A.~Rodr\'{\i}guez-Palacios}, Non associative $C^*$-algebras revisited,
In: \emph{Recent Progress in Functional analysis}, Proceedings of
the International Function Analysis Meeting on the Occasion of the
70th Birthday of Professor Manuel Valdivia (K.~D.~Bierstedt,
J.~Bonet, M.~Maestre and J.~Schmets Eds.). Elsevier, Amsterdam,
2001, pp.\ 379--408.

\bibitem{Lac} \textsc{H.~E.~Lacey},
\emph{The isometric theory of classical Banach spaces},
Springer-Verlag, Berlin, 1972.

\bibitem{Lima} \textsc{\AA.~Lima},
Intersection properties of balls and subspaces in Banach spaces,
\emph{Trans. Amer. Math. Soc.}, \textbf{227} (1977), 1--62.

\bibitem{Lima2} \textsc{\AA.~Lima},
Intersection properties of balls in spaces of compact operators,
\emph{Ann. Inst. Fourier Grenoble}, \textbf{28} (1978), 35--65.

\bibitem{Lin} \textsc{J.~Lindenstrauss},
\emph{Extension of compact operators}, Memoirs of the Amer. Math.
Soc. \textbf{48}, Providence, 1964.

\bibitem{Lin0} \textsc{J.~Lindenstrauss},
On operators which attain their norm, \emph{Israel J. Math.}
\textbf{1} (1968), 139--148.

\bibitem{L-M-M} \textsc{G.~L\'{o}pez, M.~Mart\'{\i}n, and J.~Mer\'{\i}},
Numerical index of Banach spaces of weakly or weakly-star continuous
functions, \emph{Rocky Mount. J. Math.} (to appear).

\bibitem{L-M-P} \textsc{G.~L\'{o}pez, M.~Mart\'{\i}n, and R.~Pay\'{a}},
Real Banach spaces with numerical index 1, \emph{Bull. London Math.
Soc.} \textbf{31} (1999), 207--212.

\bibitem{Lumer} \textsc{G.~Lumer},
Semi-inner-product spaces, \emph{Trans. Amer. Math. Soc.}
\textbf{100} (1961), 29--43.

\bibitem{Mar} \textsc{M.~Mart\'{\i}n},
A survey on the numerical index of a Banach space, III Congress on
Banach Spaces (Jarandilla de la Vera, 1998), \emph{Extracta Math.}
\textbf{15} (2000), 265--276.

\bibitem{MarEMA} \textsc{M.~Mart\'{\i}n},
\'{I}ndice num\'{e}rico y subespacios (Spanish), \emph{Proceedings of the
Meeting of Andalusian Mathematicians (Sevilla, 2000)}, Vol.~II, pp.\
641--648, Colecc. Abierta \textbf{52}, Univ. Sevilla Secr. Publ.,
Sevilla, 2001.

\bibitem{MarRNP} \textsc{M.~Mart\'{\i}n},
Banach spaces having the Radon-Nikod\'{y}m property and numerical index
$1$, \emph{Proc. Amer. Math. Soc.} \textbf{131} (2003), 3407--3410.

\bibitem{Mar-ADP} \textsc{M.~Mart\'{\i}n},
The alternative Daugavet property for $C^*$-algebras and
$JB^*$-triples, \emph{Math. Nachr.} (to appear).

\bibitem{MarMer} \textsc{M.~Mart\'{\i}n and J.~Mer\'{\i}}, Numerical index
of some polyhedral norms on the plane, \emph{Linear Multilinear
Algebra} (to appear).

\bibitem{MaMeRo} \textsc{M.~Mart\'{\i}n, J.~Mer\'{\i}, and A.~Rodr\'{\i}guez-Palacios},
Finite-dimensional Banach spaces with numerical index zero,
\emph{Indiana University Math. J.} \textbf{53} (2004), 1279--1289.

\bibitem{MaOi} \textsc{M.~Mart\'{\i}n and T.~Oikhberg},
An alternative Daugavet property, \emph{J. Math. Anal. Appl.}
\textbf{294} (2004), 158--180.

\bibitem{M-P} \textsc{M.~Mart\'{\i}n and R.~Pay\'{a}},
Numerical index of vector-valued function spaces, \emph{Studia
Math.} \textbf{142} (2000), 269--280.

\bibitem{MartPaya-CL} \textsc{M.~Mart\'{\i}n and R.~Pay\'{a}},
On CL-spaces and almost-CL-spaces, \emph{Ark. Mat.} \textbf{42}
(2004), 107--118.

\bibitem{M-V} \textsc{M.~Mart\'{\i}n and A.~R.~Villena},
Numerical index and Daugavet property for $L_\infty(\mu,X)$,
\emph{Proc. Edinburgh Math. Soc.} \textbf{46} (2003), 415--420.

\bibitem{Mc} \textsc{C.~M.~McGregor},
Finite dimensional normed linear spaces with numerical index $1$,
\emph{J. London Math. Soc.} \textbf{3} (1971), 717--721.

\bibitem{MPRY} \textsc{J.~F.~Mena, R.~Pay\'{a}, A.~Rodr\'{\i}guez-Palacios, and
D.~Yost}, Absolutely proximinal subspaces of Banach spaces, \emph{J.
Aprox. Theory} \textbf{65} (1991), 46--72.

\bibitem{Neg} \textsc{S.~Negrepontis}, Banach spaces and
topology, in: \emph{Handbook of set theoretic topology} (K.~Kunen
and J.~E.~Vaughan, eds.), North Holland, Amsterdam, 1984.

\bibitem{Oik04} \textsc{T.~Oikhberg},
Spaces of operators, the $\psi$-Daugavet property, and numerical
indices, \emph{Positivity} \textbf{9} (2005), 607--623.

\bibitem{Par} \textsc{J.~R.~Partington,}
Norm attaining operators, \emph{Israel J. Math.} \textbf{43} (1982),
273--276.

\bibitem{Pelcz} \textsc{A.~Pe\l czy\'nski},
On the impossibility of embedding of the space $L$ in certain Banach
spaces, \emph{Colloq. Math.} \textbf{8} (1961), 199--203.

\bibitem{Reis} \textsc{S.~Reisner},
Certain Banach spaces associated with graphs and CL-spaces with
$1$-unconditonal bases, \emph{J. London Math. Soc} \textbf{43}
(1991), 137--148.

\bibitem{Rolewicz} \textsc{S.~Rolewicz},
\emph{Metric Linear Spaces (Second edition)}, Mathematics and its
Applications (East European Series), \textbf{20}. D. Reidel
Publishing Co., Dordrecht; PWN---Polish Scientific Publishers,
Warsaw, 1985.

\bibitem{Rosenthal} \textsc{H.~P.~Rosenthal}, The Lie algebra of a Banach space,
in: \emph{Banach spaces} (Columbia, Mo., 1984), 129--157, Lecture
Notes in Math. \textbf{1166}, Springer, Berlin, 1985.

\bibitem{Ros-Pacific} \textsc{H.~P.~Rosenthal}, Functional hilbertian
sums, \emph{Pac. J. Math.} \textbf{124} (1986), 417--467.

\bibitem{Scha} \textsc{W.~Schachermayer},
Norm attaining operators and renormings of Banach spaces,
\emph{Israel J. Math.} \textbf{44} (1983), 201--212.

\bibitem{Sch} \textsc{K.~D.~Schmidt},
Daugavet's equation and orthomorphisms, \emph{Proc. Amer. Math.
Soc.} \textbf{108} (1990), 905--911.

\bibitem{Shv} \textsc{R.~V.~Shvidkoy},
Geometric aspects of the Daugavet property, \emph{J. Funct. Anal.}
\textbf{176} (2000), 198--212.

\bibitem{Tille} \textsc{K.~Tillekeratne},
Spatial numerical range of an operator, \emph{Proc. Camb. Phil.
Soc.} \textbf{76} (1974), 515--520.

\bibitem{Toe} \textsc{O.~Toeplitz},
Das algebraische Analogon zu einem Satze von Fejer, \emph{Math. Z.}
\textbf{2} (1918), 187--197.

\bibitem{Val} \textsc{M.~Valdivia},
Topological direct sum decompositions of Banach spaces, \emph{Israel
J. Math.} \textbf{71} (1990), 289--296.

\bibitem{Wer0} \textsc{D.~Werner},
An elementary approach to the Daugavet equation, in:
\emph{Interaction between Functional Analysis, Harmonic Analysis and
Probability} (N.~Kalton, E.~Saab and S.~Montgomery-Smith Eds),
Lecture Notes in Pure and Appl. Math. \textbf{175} (1994), 449--454.

\bibitem{WerJFA} \textsc{D.~Werner},
The Daugavet equation for operators on function spaces, \emph{J.
Funct. Anal.} \textbf{143} (1997), 117--128.

\bibitem{WerSur} \textsc{D.~Werner},
Recent progress on the Daugavet property, \emph{Irish Math. Soc.
Bull.} \textbf{46} (2001), 77--97.

\bibitem{Woj} \textsc{P.~Wojtaszczyk},
Some remarks on the Daugavet equation, \emph{Proc. Amer. Math. Soc.}
\textbf{115} (1992), 1047--1052.


\end{thebibliography}
\end{document}